\documentclass[11pt,graphicx]{article}
\usepackage{amsmath}
\usepackage{amsfonts}
\usepackage{color}

\def\X{{\cal X}}

\def\T{{\cal T}}

\newfont{\Blackboard}{msbm10 scaled 1200}

\newfont{\roma}{cmr10 scaled 1200}

\newtheorem{thm}{{}\hskip\parindent Theorem}[section]
\newtheorem{lem}{{}\hskip\parindent Lemma}[section]
\newtheorem{pro}{{}\hskip\parindent Proposition}[section]
\newtheorem{exl}{{}\hskip\parindent Example}[section]
\newtheorem{cor}{{}\hskip\parindent Corollary}[section]
\newtheorem{definition}{{}\hskip\parindent Definition}[section]
\newtheorem{rem}{{}\hskip\parindent Remark}[section]

\def\be{\begin{equation}}
\def\ee{\end{equation}}
\def\beq{\arraycolsep=1.5pt\begin{eqnarray}}
\def\eeq{\end{eqnarray}}

\large
\title{Second order necessary  conditions for optimal control problems  with endpoints-constraints and convex control-constraints \thanks{This work is  supported by the National Science
Foundation of China under grants 11401491 and 11231007, and the Fundamental research funds for the Central Universities under grant 2682014CX052.
}
\date{}
\author{Li Deng\thanks{School of Mathematics,   Southwest Jiaotong University, Chengdu 611756, Sichuan Province, China. {\small\it E-mail:} {\small\tt
dengli@swjtu.edu.cn}.}}
}

\begin{document}
\maketitle

\begin{quote}
\begin{small}
{\bf Abstract} \,\,\,
In this manuscript, we consider a control system governed by a general ordinary differential equation on a Riemannian manifold, with its endpoints satisfying some inequalities and equalities, and its control constrained to a closed convex set. We concern
on an optimal control problem of this system, and obtain the second order necessary
condition in the sense of convex variation (Theorem \ref{c230}). To this end, we first obtain
a second order necessary condition of an optimization problem (Theorem \ref{c126}) via  separation theorem of convex sets. Then, we derive our necessary condition by transforming
the optimal control problem into an optimization problem.
It is worth to point out that, our necessary condtition evolves the curvature tensor,
which is trivial in Euclidean case. Moreover, even M is a Euclidean space, our result is
still of interest. Actually, we give an example (Example \ref{ccs126}) which shows that, when
an optimal control stays at the boundary of the control set, the existing results are
invalid while Theorem \ref{c230} works.
\\[3mm]
{\bf Keywords}\,\,\, Optimal control, Second order necessary   condition, Endpoints-constraints, Convex constraints, Riemannian manifold
\\[3mm]
{\bf MSC (2010) \,\,\,
49K15, 49K30, 93C15, 58E25, 70Q05} \\[3mm]
\end{small}
\end{quote}

\setcounter{equation}{0}

\section{Introduction }
\def\theequation{1.\arabic{equation}}
\hskip\parindent

In this paper, we consider a control system described by a general ordinary differential equation with the state restricted to a manifold, with the initial and terminal states restricted to inequality-type and equality-type constraints, and with the control constrained pointwisely to a convex set. For this control system, we study the second order necessary optimality condition for an optimal control problem.

Before elaborating on our problem, we introduce some notions on manifolds.
 Let $n\in\mathbb N$ and $M$ be a complete simply connected, $n$-dimensional manifold with Riemannian metric $g$. Let $\nabla$ be the Levi-Civita connection on $M$ related to $g$, $\rho(\cdot,\cdot)$ be the distance function on $M$, and  $T_xM$  and $T^*_x M$ be respectively the tangent and cotangent spaces of  $M$ at $x\in M$. Denote by $\langle\cdot,\cdot\rangle$ and $|\cdot|$ the inner product and the norm over $T_xM$ related to $g$, respectively. Also, denote by  $ T M \equiv \bigcup\limits_{x\in M} T_xM$, $T^* M \equiv \bigcup\limits_{x\in M} T^*_x M$ and $C^\infty(M)$ the tangent bundle, the cotangent bundle and the set of smooth functions on $M$, respectively.

Let $j,k\in\mathbb N$, $T>0$, $U$ be a subset of $\mathbb R^m$ ($m\in\mathbb N$), and  $f:[0,T]\times M\times U\to T M$, $\phi_i:M\times M\to \mathbb R$ ($i=0,1,\cdots,j$) and $\psi=(\psi_1, \cdots, \psi_k)^\top:M\times M\to \mathbb R^k$ be maps (satisfying suitable assumptions to be given later).   Set by
\begin{align}\ \label{cs}& \mathcal U=\{v:[0,T]\to U|\,v(\cdot)\,\textrm{is measurable}\}
\end{align}
the set of all possible controls.
Consider the following optimal control problem 
\begin{description}
\item[$(OCP)$] Find  a control $\bar u(\cdot)$ and a trajectory $\bar y(\cdot)$   minimizing
$$
J(y(\cdot), u(\cdot))\equiv \phi_0(y(0),y(T)),
$$
subject to
\begin{equation}\label{s2}
\dot{y}(t)=f(t,y(t),u(t)),\;a.e.\,t\in[0,T];\;u(\cdot)\in\mathcal U;
\end{equation}
and 
\begin{equation}\label{s1}
\left\{\begin{array}{l}
 \phi_i(y(0),y(T))\leq 0,\; i=1,\cdots,j,
\\ \psi(y(0),y(T))=0.
\end{array}\right.
\end{equation}
$\bar u(\cdot)$, $\bar y(\cdot)$ and $(\bar u(\cdot),\bar y(\cdot))$ are respectively called  optimal control, optimal trajectory and optimal pair.
\end{description}

For problem $(OCP)$,
  Deng and Zhang (\cite{dzgJ})   obtained the second order necessary  condition via  spike variation, by applying the separation theorem of convex sets to a suitably chosen set related to the second order spike variation. In this paper,   we are concerned with the second order necessary condition obtained by convex variation.

In optimal control theory, optimality conditions are usually obtained by two kinds of variations of trjectories: spike  and convex variations. In this paper, we call a necessary condition obained by convex variation (resp. spike variation) a necessary condition in the sense of convex variation (resp. spike variation).
For the first order necessary optimality condition, the condition  in the sense of spike variation is always more precise than that in the sense of convex variation, i.e. the latter can be deduced directly from the former. However, for the second order necessary condition, it is hard to assess which condition is better. As already explained in \cite[Section 1]{dzgJ}, the second order necessary condition only makes sense along a critical direction, in which   the first order necessary condition is trivial. Since different variations lead to different first order necessary conditions, critical directions in the senses of different variations are distinct. Consequently the corresponding second order necessary conditions are different. They are always complementary to each other. 
Two examples in \cite{War78} (or \cite[Section 3.5]{dzs}) show that, these two kinds of second order necessary conditions can not cover each other. Our main result ( Theorem \ref{c230}) is the second order necessary condition  in the sense of convex variation, which can be viewed as a complement to \cite[Theorem 2.2]{dzgJ}.

For problem $(OCP)$, when $M$ is a Euclidean space, the second order necessary condition in the sense of convex variation is given by \cite[Theorem 3.1]{gb}. Comparing to it, our main result (i.e. Theorem \ref{c230}) has two differences: i)  The curvature tensor of the Riemannian manifold arises. This is the same case as the corresponding second order necessary condition in the sense of spike variation (i.e. \cite[Theorem 2.2]{dzgJ}); ii) There is another extra term  (i.e. the first integral of the left hand side of inequality (\ref{c244})), which always makes sense when an optimal control stays at the boundary of the control set, and provides additional informations. We explain this by Example \ref{ccs126}, in which Theorem \ref{c230} works, while \cite[Theorem 3.1]{gb} fails.

The appearance of the curvature tensor results from the second order variation of trajectories evolved on Riemannian manifolds.  The extra term (the first integral of the left hand side of inequality (\ref{c244})) comes from the use of ``the second-order adjacent subset" introduced in  \cite[Definition 4.7.2, p. 171]{AF}. More precisely,
we prove Theorem \ref{c230} by using the second order necessary condition of a solution to an optimization problem with inequality-type and equality-type constraints (i.e. Theorem \ref{c126}), which evolves ``the second-order adjacent subset". With this nontrivial term, Theorem \ref{c126}  generalizes \cite[Theorem 4.1]{gb}, see Remark \ref{c225} for details.

This paper is organised as follows. The main results are stated in Section 2, the variations of trajectories to the second order are given in Section 3, and Section 4 is devoted to the proof of the main results.

\setcounter{equation}{0}
\hskip\parindent \section{Statement of the main results}
\def\theequation{2.\arabic{equation}}

\subsection{Notations and assumptions}\label{na}
We first introduce some notions. Denote by $i(x)$, $|\T(x)|$, $\nabla \T$, $R$, the injectivity radius (at the point  $x\in M$), the norm of the tensor field $\T$ at the point $x\in M$ (see \cite[(2.5)]{cdz}), the covariant derivative of the tensor field $\T$ and the curvature tensor ( of $(M,g)$), respectively. For any $x,y\in M$ with $\rho(x,y)<\min\{i(x),i(y)\}$, by \cite[Lemma 2.1]{cdz}, there exists a unique shortest geodesic connecting $x$ and $y$. We denote   the parallel translation of a tensor from $x$ to $y$ along this geodesic  by   $L_{xy}$.  For a smooth function $h: M\times M\to \mathbb R$ of two arguments, we denote by $\nabla_i h(y_1,y_2)$  the covariant derivative of $h$ with respect to  the $i^{th}$ argument $y_i$ with $i=1,2$. Namely, we have
\begin{align*}
\nabla_ih(y_1,y_2)(X)=X(y_i)(h(y_1,y_2)),\;\forall\,X\in T M.
\end{align*}
Thus, $\nabla_ih(y_1,y_2)\in T_{y_i}^*M$.  For a smooth vector valued map $\psi=(\psi_1,\cdots,\psi_k)^\top:M\times M\to \mathbb R^k$ ($k>0$), we denote  the first and   second order covariant derivatives with respect to the $i^{th}$ argument $(i=1,2)$ respectively by
\begin{align}\label{pcd}\begin{array}{l}
\nabla_i\psi(X)=(\nabla_i\psi_1(X),\cdots,\nabla_i\psi_k(X))^\top,
\\[2mm]\nabla_i^2\psi(X,Y)=(\nabla_i^2\psi_1(X,Y),\cdots,\nabla_i^2\psi_k(X,Y))^\top, \end{array}
\end{align}
for all $X,Y\in T M.$
The corresponding norms are given respectively by $|\nabla_i\psi|=\sum_{l=1}^k|\nabla_i\psi_l|$ and $|\nabla_i^2\psi|=\sum_{l=1}^k|\nabla_i^2\psi_l|$. For the definitions of the above notions, please see \cite[Section. 2]{cdz}.

To present the optimality conditions of problem $(OCP)$, we need to introduce two functions: Lagrange and  Hamiltonion functions. When $k>0$, the Lagrange function  $\mathcal L:M\times M\times \mathbb R^{1+j+k}\to \mathbb R$  is defined by $\mathcal L(y_1,y_2;\ell)\equiv\sum_{i=0}^j\ell_i\phi_i(y_1,y_2)+\ell_\psi^\top\psi(y_1,y_2)$, where $\ell=(\ell_0,\cdots,\ell_j,\ell_\psi^\top)^\top$. When $k=0$, the corresponding Lagrange function  $\mathcal L:M\times M\times \mathbb R^{1+j}\to \mathbb R$ is $\mathcal L(y_1,y_2;\ell)\equiv\sum_{i=0}^j\ell_i\phi_i(y_1,y_2)$, where $\ell=(\ell_0,\cdots,\ell_j)^\top$.  For each $\ell\in \mathbb R^{1+j+k}$,
we also denote by $d_i\mathcal L(y_1,y_2;\ell)$ the differential of $\mathcal L$ with respect to the variable $y_i$ ($i=1,2$), i.e. $d_i\mathcal L(y_1,y_2;\ell) \in T_{y_i}^*M$  satisfies $d_i\mathcal L(y_1,y_2;\ell) (X)= X(y_i)(\mathcal L(y_1,y_2;\ell))$ for all $X\in T M$.
The Hamiltonian function $H:[0,T]\times T^*M\times U\to\mathbb R$ is given by 
\begin{equation}\label{c250}H(t,y,p,u)\equiv p(f(t,y,u)),\quad\forall\,(t,y,p,u)\in[0,T]\times T^*M\times U.\end{equation}

Throughout this paper, we denote by  $\tilde X\in T^*M$ (resp. $\tilde X\in T M$) the dual covector (resp. vector) of  $X\in T M$ (resp. $X\in T^*M$), see \cite[Section 2.1]{dzgJ} for the detailed definition.

Then, we recall some definitions concerning tangent sets. For more details, we refer to \cite{AF}. Let  $\mathcal X$ be a metric space with a metric $d$, and $K\subset \mathcal X$ be a subset. The distance between a point $x\in\mathcal X$ and $K$ is defined by $dist_K(x):=\inf\{d(x,y); y\in K\}$. Let $\{K_h\}_{h>0}$ be a family of subsets of $\mathcal X$. The lower limit of $\{K_h\}_{h>0}$ is given by 
$$
{Liminf}_{h\to 0^+}K_h:=\{v\in \mathcal X;\,\lim_{h\to 0^+} dist_{K_h}(v)=0\}.
$$
When $\mathcal X$ is a normed vector space,
the adjacent cone to a subset $K\subset \mathcal X$ at a point $x\in \overline K$ (i.e. $x$ belongs to the closure of $K$) is defined by (see \cite[p.127]{AF})
\begin{align}\label{ccs94}
T_K^\flat(x):=Liminf_{h\to 0^+}\frac{K-x}{h}.
\end{align}
Moreover, for $v\in T_K^\flat(x)$, the second-order adjacent subset to $K$ at $(x,v)$ is given by (see \cite[Definition 4.7.2, p. 171]{AF})
\begin{align}\label{ccs90}
T_K^{\flat(2)}(x,v):=Liminf_{h\to 0^+}\frac{K-x-hv}{h^2}.
\end{align}
By the definition of the lower limit of a family of subsets, one can  respectively characterise $T_K^\flat(x)$ and $T_K^{\flat(2)}(x,v)$ in terms of sequences:
\begin{description}
\item[(i)] $v\in T_K^\flat(x)$ if and only if, for any $h_n\to 0^+$ as $n\to+\infty$, there exists $v_n\in \mathcal X$ approaching to $v$, such that $x+h_n v_n\in K$ for each $n\geq 1$;
\item[(ii)] Given $v\in T_K^\flat(x)$, $w\in T_K^{\flat(2)}(x,v)$ if and only if,  for any $h_n\to 0^+$ as $n\to+\infty$, there exists $w_n\in\mathcal X$ approaching to $w$, such that $x+h_nv+h_n^2w_n\in K$ for each $n\geq 1$.
\end{description}

The main assumptions are exhibited as follows:
\begin{description}
\item[$(C1)$] $U\subset \mathbb R^m$ is convex and closed.
\item[$(C2)$] The map $f(=f(t,x,u)): [0, T]\times M\times U \to T M$ is measurable in $t$,  and $C^1$ in $(x,u)$. Moreover,   there exists a constant $K>1$ such that,
\begin{equation}\begin{array}{l}\label{10}
 |L_{x_1\hat x_1}f(s,x_1,u_1)-f(s,\hat x_1,u_2)|\leq K(\rho(x_1,\hat x_1)+|u_1-u_2|),
\\|f(s,x_0,u_1)|\leq K,
\\ |\phi_i(x_1,x_2)-\phi_i(\hat x_1,\hat x_2)|\leq K(\rho(x_1,\hat x_1)+\rho(x_2,\hat x_2)),\;i=0,\cdots,j,
\\|\psi(x_1,x_2)-\psi(\hat x_1,\hat x_2)|\leq K (\rho(x_1,\hat x_1)+\rho(x_2,\hat x_2)),
\end{array}
\end{equation} for all $s\in [0,T]$, $u_1,u_2\in U$, and  $x_l, \hat x_l\in M$ satisfying  $\rho(x_l,\hat x_l)< \min\{i(x_l),i(\hat x_l)\}$ for $l=1,2$,  where $x_0\in M$ is  fixed.

\item[$(C3)$]The map $f=(f(s,x,u))$ is $C^2$ in $(x,u)\in M\times U$, and $\phi_0, \cdots,\phi_j, \psi$ are $C^2 $ over $M\times M$.   Furthermore,  there exists a positive constant $K>1$ such that
 \begin{equation}\label{116}\begin{array}{l}
|\nabla_{x}f(s,x_1,u_1)-L_{\hat x_1x_1}\nabla_x f(s,\hat x_1,u_2)|\leq K(\rho(x_1,\hat x_1)+|u_1-u_2|),
\\ |\nabla_uf(s,x_1,u_1)-L_{\hat x_1x_1}\nabla_uf(s,\hat x_1,u_2)|\leq K(\rho(x_1,\hat x_1)+|u_1-u_2|),
\\| \nabla_1\phi_i(x_1,x_2)-L_{\hat x_1 x_1}\nabla_1\phi_i(\hat x_1, x_2)|\leq K \rho(x_1,\hat x_1),\;i=0, 1, \cdots, j,
\\ | \nabla_2\phi_i(x_1,x_2)-L_{\hat x_2 x_2}\nabla_2\phi_i( x_1, \hat x_2)|\leq K \rho(x_2,\hat x_2),\;i=0, 1, \cdots, j,
\\| \nabla_1\psi(x_1,x_2)-L_{\hat x_1 x_1}\nabla_1\psi(\hat x_1, x_2)|\leq K \rho(x_1,\hat x_1),
\\ | \nabla_2\psi(x_1,x_2)-L_{\hat x_2 x_2}\nabla_2\psi( x_1, \hat x_2)|\leq K \rho(x_2,\hat x_2),
\end{array}\end{equation} for all $x_1,\hat x_1, x_2,\hat x_2\in M$ with $\rho(x_1,\hat x_1)<\min\{i(x_1),i(\hat x_1)\}$ and $\rho(x_2,\hat x_2)<\min\{i(x_2),
\\ i(\hat x_2)\}$, and $(s,u_1,u_2)\in[0,T]\times U\times U$,   where 
 $\nabla_uf(s,x,u)$ is defined by 
\begin{align}\label{ccs122}\begin{array}{ll}
& \nabla_uf(s,x,u)(\eta,V)
\\ =&\lim_{\epsilon\to 0^+}
\frac{1}{\epsilon}\Big(f(s,x,u+\epsilon V)(\eta)-f(s,x,u)(\eta)\Big),\;\forall\, (\eta, V)\in T_x^*M\times\mathbb R^m,
\end{array}
\end{align}
with its norm
\begin{align*}
|\nabla_uf(s,x,u)|\equiv\sup\{\nabla_uf(s,x,u)(\eta,V);\,(\eta,V)\in T_x^*M\times\mathbb R^m, |\eta|+|V|\leq 1\},
\end{align*}
 and $\nabla_x f(s,x,u)$ is  the covariant derivative of  $f(s,x,u)$    with respect to the state  variable $x\in M$, and  is    defined by
\begin{align}\label{ccs120}
\nabla_x f(s,x,u)(\eta,X)=\nabla_X f(s,\cdot,u)(\eta),\quad\forall\,(\eta,X)\in T_x^*M\times T_xM,
\end{align}
with its norm  given by \cite[(2.5)]{cdz}.

\end{description}

We should mention that, for  $l=1,2$ and $ i=0,\cdots,j$, by \cite[Lemma 4.1]{cdz}, we obtain from   (\ref{10}) and (\ref{116}) that,  $f$ and $\nabla_x f$ are both Lipschitz continuous with respect to the variable $(x,u)\in M\times U$, and $\phi_i$, $\psi$,  $\nabla_l\phi_i$  and $\nabla_l\psi$  are Lipschitz continuous. These conditions  can be checked by computing the norms of $\nabla_xf$, $\nabla_x^2f$, $\nabla_l\phi_i$, $\nabla_l^2\phi_i$, $\nabla_l\psi$ and $\nabla_l^2\psi$.

\subsection{Main results}\label{mr}

In this subsection, we fix an optimal pair $(\bar u(\cdot),\bar y(\cdot))$.   For a function defined on $[0,T]\times M\times U$, 
 we denote by
\begin{equation}\label{434}
\varphi [t]\equiv\varphi (t,\bar y(t),\bar u(t)),\quad\forall\, t\in[0,T]
\end{equation}
for abbreviation. 
Set 
\begin{eqnarray}\label{c246}
I_{A}\equiv\{i\in\{1,\cdots,j\};\,\phi_i(\bar y(0),\bar y(T))=0\}\cup\{0\},
\\ \label{c226}
I_N\equiv\{0,1,\cdots,j\}\setminus I_{A}.
\end{eqnarray}
Given a vector $\ell=(\ell_0,\cdots,\ell_j,\ell_\psi^\top)^\top\in \mathbb R^{1+j+k}$, we denote by $p^\ell(\cdot)$ the solution to 
\begin{equation}\label{c243}
\left\{\begin{array}{l}
\nabla_{\dot{\bar y}(t)}p^\ell=-\nabla_xf[t](p^\ell(t),\cdot),\;a.e.\,t\in(0,T),
\\ p^\ell(T)=d_2\mathcal L(\bar y(0),\bar y(T);\ell),
\end{array}\right.
\end{equation}
with $\nabla_xf$ given  by (\ref{ccs120}). It is a covector field along $\bar y(\cdot)$, i.e. $p^\ell(t)\in T^*_{\bar y (t)} M$ for each $t\in[0,T]$. Furthermore, for a function $\varphi$ defined on $[0,T]\times T^*M\times U$, we  set
\begin{equation}\label{c245}
\varphi[t,\ell]=\varphi (t,\bar y(t),p^\ell(t),\bar u(t)),\;\forall\, t\in[0,T]
\end{equation}
for abbreviation.

The first order necessary condition of an optimal pair in the secnse of convex variation is stated as follows.
\begin{thm}\label{f}
Assume $U\subset\mathbb R^m$ ($m\in\mathbb N$) is convex, and condition $(C2)$  holds. If $( \bar u(\cdot),\bar y(\cdot))$ is optimal pair for problem  $(OCP)$ with $\bar u(\cdot)\in L^2(0,T;\mathbb R^m)\cap \mathcal U$, then there exists $\ell=(\ell_0,\ell_1,\cdots,\ell_j,\ell_\psi^\top)^\top\in\mathbb R^{1+j+k}\setminus\{0\}
$
such that
\begin{equation}\label{c260}\begin{array}{l}
\ell_i\leq 0,\; \textrm{if}\;i\in I_A;\quad \ell_i=0,\;\textrm{if}\;i\in I_N,
\end{array}\end{equation}
and 
\begin{equation}\label{c240}
\begin{array}{l}
\nabla_u H[t,\ell](v(t))\leq 0,\;a.e. \,t\in[0,T],
\end{array}
\end{equation}
holds for all $v(\cdot)\in L^2(0,T;\mathbb R^m)$ with $v(t)\in T_U^\flat(\bar u(t))$ a.e $t\in[0,T]$,
where $p^\ell(\cdot)$ is a covector field along $\bar y(\cdot)$ satisfying (\ref{c243}) and initial condition
\begin{align}\label{icd1}
p^\ell(0)=-d_1\mathcal L(\bar y(0),\bar y(T);\ell),
\end{align}
and $\nabla_u H(t,y,\eta,u)(V)$ (with $(t,y,\eta,u,V)\in[0,T]\times T^*M\times U\times\mathbb R^m$) is defined by
\begin{align}\label{ccs124}
\nabla_uH(t,y,\eta,u)(V)=\lim_{\epsilon\to0^+}\frac{1}{\epsilon}\Big(H(t,y,\eta,u+\epsilon V)-H(t,y,\eta,u)\Big).
\end{align}
\end{thm}

From the viewpoint of calculus, when the first order necessary condition is trivial in some direction, it is necessary to find the second order necessary condition along this direction. Thus, in what follows, we give the definition of critical direction in the sense of convex variation.

\begin{definition}\label{c228} Assume that $( \bar u(\cdot), \bar y(\cdot))$ is an optimal pair of problem $(OCP)$ with $\bar u(\cdot)\in L^2(0,T;\mathbb R^m)\cap\mathcal U$, and that all the assumptions in Theorem \ref{f} hold.
A function $v(\cdot)\in  L^2(0,T;\mathbb R^m)$ is called a singular direction in the sense of convex variation, if it satisfies
\begin{equation}\label{c266}\begin{array}{l}
v(t)\in T_U^\flat(\bar u(t)),\quad a.e.\,t\in[0,T],
\\[2mm]
\nabla_1\phi_i(\bar y(0),\bar y(T))(X_v(0))+\nabla_2\phi_i(\bar y(0),\bar y(T))(X_{v}(T))\leq 0,\;\forall \,i\in  I_{A},
\\[2mm]
\nabla_1\psi(\bar y(0),\bar y(T))(X_v(0))+\nabla_2\psi(\bar y(0),\bar y(T))(X_{v}(T))=0\quad(\textrm{omit if}\;k=0),
\end{array}
\end{equation}
where $\nabla_i\psi$ is defined by (\ref{pcd}),
and $X_{v}(\cdot)$ is a vector field along $\bar y(\cdot)$ (i.e. $X_v(t)\in T_{\bar y(t)}M$ for all $t\in[0,T]$) and verifies
 \begin{equation}\label{c227}
\nabla_{\dot{\bar y}(t)}X_{v}=\nabla_xf[t](\cdot,X_{v}(t))+\nabla_uf[t](\cdot,v(t)),\;a.e.\,t\in(0,T),
\end{equation}
with $\nabla_xf$ and $\nabla_uf$ given respectively  by (\ref{ccs120}) and (\ref{ccs122}).

Along a critical direction $v(\cdot)$ defined above, the first order necessary condition in the sense of convex variation is trival. To show this, we need a definition as follows.

\begin{definition}\label{lpc}
Assume that all the assumptions in Theorem \ref{f} hold, and that $(\bar u(\cdot),\\ \bar y(\cdot))$ is an optimal pair of problem $(OCP)$. A vector $\ell=(\ell_0,\ell_1,\cdots,\ell_j,\ell_\psi^\top)^\top\in\mathbb R^{1+j+k}\setminus\{0\}$ is called a Lagrange multiplier  in the sense of convex variation, if it satisfies (\ref{c260}), and  (\ref{c240}) holds for all $v(\cdot)\in L^2(0,T;\mathbb R^m)$ with $v(t)\in T_U^\flat(\bar u(t))$ a.e $t\in[0,T]$,
where $p^\ell(\cdot)$ satisfies (\ref{c243}) and  (\ref{icd1}).
\end{definition}

Thus, we can understand
the first order necessary condition  in the sense of convex variation (i.e. Theorem \ref{f}) as follows: if $(\bar u(\cdot),\bar y(\cdot))$ is an optimal pair, there exists a Lagrange multiplier  in the sense of convex variation. Moreover, if $v(\cdot)\in L^2(0,T;\mathbb R^m)$ is a crictical direction in the sense of convex variation, with $X_v(\cdot)$ satisfying (\ref{c266}) and (\ref{c227}),  then for any Lagrange multiplier $\ell\in\mathbb R^{1+j+k}\setminus\{0\}$, by (\ref{c266}),  (\ref{c243}), (\ref{icd1}), (\ref{c227}), (\ref{c240}) and integration by parts, we have
\begin{align*}
0\leq &\nabla_1\mathcal L(\bar y(0),\bar y(T);\ell)(X_v(0))+\nabla_2\mathcal L(\bar y(0),\bar y(T);\ell)(X_v(T))
\\=&\int_0^T\nabla_uH[t,\ell](v(t))dt\leq0,
\end{align*}
which implies $\nabla_uH[t,\ell](v(t))=0$ a.e $t\in[0,T]$.
Thus, along direction $v(\cdot)$, the first order necessary condition in the sense of convex variation  is trivial.

Then, along this direction $v(\cdot)$, we shall study the second order necessary condition. To this end,
associated to $v(\cdot)$ and $X_v(\cdot)$, we set
\begin{equation}\label{c231}\begin{array}{l}
 I_0^\prime\equiv  I_N\cup\{i\in I_{A};\,\nabla_1\phi_i(\bar y(0),\bar y(T))(X_v(0))+\nabla_2\phi_i(\bar y(0),\bar y(T))(X_{v}(T))<0\};
\\[2mm]  I_0^{\prime\prime}\equiv \{0,1,\cdots,j\}\setminus  I_0^\prime.
\end{array}
\end{equation}
\end{definition}
The second order necessary condition of optimal pairs is as follows.
\begin{thm}\label{c230}
Assume that conditions $(C1)-(C3)$ hold, and that $(\bar u(\cdot), \bar y(\cdot))$ with $\bar u(\cdot)\in L^2(0,T;\mathbb R^m)\cap\mathcal U$ is an optimal pair of problem $(OCP)$.  Let $v(\cdot)\in   L^2(0,T;\mathbb R^m)$ be a critical  direction in the sense of convex variation, with $X_v(\cdot)$ satisfying (\ref{c266}) and (\ref{c227}). 
Assume that there exist $\ell(\cdot)\in L^2(0,T;\mathbb R^m)$ and $\epsilon_0>0$ such that 
\begin{equation}\label{l4}dist_U(\bar u(t)+\epsilon v(t))\leq  \epsilon^2\ell(t),\forall \epsilon\in[0,\epsilon_0], a.e.\,t\in[0,T],\end{equation} 
and that the set $\mathcal B\equiv\{\sigma(\cdot)\in L^2(0,T;\mathbb R^m);\,\sigma(t)\in T_U^{\flat(2)}(\bar u(t),v(t))\,a.e.\,t\in[0,T]\}\not=\emptyset$. 
Then, there exists a Lagrange multiplier $\ell=(\ell_{\phi_0},\ell_{\phi_1},\cdots,\ell_{\phi_j},\ell_\psi^\top)^\top\in\mathbb R^{1+j+k}\setminus\{0\}$ satisfying 
\begin{equation}\label{c232}
 \ell_{\phi_i}=0,\;\textrm{if}\;i\notin  I_0^{\prime\prime},
\end{equation}
 such that 
\begin{equation}\label{c244}
\begin{array}{l}
\int_0^T\nabla_uH[t,\ell](\sigma (t))dt+\frac{1}{2}\int_0^T\left\{\nabla_x^2H[t,\ell](X_{v}(t),X_{v}(t))\right.
\\+2\nabla_u\nabla_xH[t,\ell](X_{v}(t),v(t))+\nabla_u^2H[t,\ell](v(t),v(t))
\\\left. -R(\tilde{ p}^\ell(t),X_{v}(t),f[t],X_{v}(t))\right\}dt
+\frac{1}{2}\nabla_1^2\mathcal L(\bar y(0),\bar y(T);\ell)(X_v(0),X_v(0))
\\+\nabla_1\nabla_2\mathcal L(\bar y(0),\bar y(T);\ell)(X_{v}(T),X_v(0))
\\+\frac{1}{2}\nabla_2^2\mathcal L(\bar y(0),\bar y(T);\ell)(X_{v}(T),X_{v}(T))\leq0
\end{array}
\end{equation}
holds for all $\sigma(\cdot)\in\mathcal B$, where $ p^{\ell}(\cdot)$  satisfies (\ref{c243}) and (\ref{icd1}), $\tilde{ p}^{\ell}(t)$ ($t\in[0,T]$) is the dual vector of $ p^{\ell}(t)$,    $H[t,\ell]$ is defined in (\ref{c250}) and (\ref{c245}),  for each $(t,y,\eta,u)\in[0,T]\times T^*M\times U$,
\begin{align*}
&\nabla_xH(t,y,\eta,u)(X)=\nabla_xf(t,y,u)(\eta,X),\quad \forall\,X\in  T_yM,
\\&\nabla_u\nabla_xH(t,y,\eta,u)(X,V)=\frac{d}{ds}\Big|_{s=0}\nabla_xH(t,y,\eta,u+sV)(X),\;\forall\,(X,V)\in  T_y M\times \mathbb R^m,
\\& \nabla_u^2H(t,y,\eta,u)(V,V)=\frac{d^2}{ds^2}\Big|_{s=0}H(t,y,\eta,u+sV),\quad\forall\,V\in\mathbb R^m,
\\&\nabla_x^2H(t,y,\eta,u)(X,X)=\nabla_x^2f(t,y,u)(\eta,X,X),\quad\forall\,X\in T_yM,
\end{align*}
and $\nabla_i\nabla_j\mathcal L(\bar y(0),\bar y(T);\ell)$ ($i,j=1,2$) is defined by \cite[(5.11)]{dzgJ}.
\end{thm}
\begin{rem}
\cite[Theorem 3.1]{gb} considers a special case of problem $(OCP)$: $M$ is a Euclidean space. 
Theorems \ref{f} $\&$ \ref{c230}  differ from \cite[Theorem 3.1]{gb} in three aspects. First, Theorem \ref{f} extends \cite[Theorem 3.1]{gb} from a Euclidean space to a Riemannian manifold. What comes new is that, the curvature tensor of the Riemannian manifold appears in the second order necessary condition. Second, \cite[Theorem 3.1]{gb}  says that, if $\bar u(\cdot)$ is an optimal control, the first order necessary condtion is that, there exists a nontrivial vector $\ell=(\ell_{\phi_0},\cdots,\ell_{\phi_j},\ell_\psi^\top)^\top\in \mathbb R^{1+j+k}\setminus\{0\}$ satisfying (\ref{c260}), such that the following inequality holds:
\begin{align}\label{ex2}
\int_0^T\nabla_u H[t,\ell](v(t))dt\leq 0,\;\forall\,v(\cdot)\in L^\infty(0,T;\mathbb R^m)\cap (\mathcal U-\{\bar u(\cdot)\}).
\end{align}
While (\ref{c240}) is of pointwise form, which is easier to be checked.  Third,  when there exists a set $A\subset[0,T]$ with its Lebesgue measure bigger than zero, such that $\bar u(t)$ belongs to the boundary of $U$ for all $t\in A$, the set $\{v(\cdot)\in L^2(0,T;\mathbb R^m)|\, v(t)\in T_U^\flat(\bar u(t))\,a.e.\,t\in[0,T]\}$ is some times larger than $L^\infty(0,T;\mathbb R^m)\cap(\mathcal U-\bar u(\cdot))$ (see Example \ref{ccs126}). Consequently, when a singular direction $v(\cdot)$ in the sense of convex variation satisfies $v(t)\in T_U^\flat(\bar u(t))\setminus\{U-\{\bar u(t)\}\}$ for  $t\in A$, compared  to \cite[Theorem 3.1 ]{gb}, we still have further information about an optimal  pair (see (\ref{c244})). We shall use Example \ref{ccs126} below to illustrate it more explicitly. 
\end{rem}

\begin{exl}\label{ccs126} Given $T\in(0,3-\sqrt 5)$ and $\theta>2$, consider the control  system
\begin{align}\label{ex1}
\left\{\begin{array}{ll}
\dot y_1(t)=u_2(t),& a.e.\,t\in(0,T),
\\ \dot y_2(t)=-y_1^2(t)+4y_1(t)u_2(t)-\theta u_1(t)^2,& a.e.\,t\in(0,T),
\end{array}
\right.
\end{align}
where $(u_1(t),u_2(t))\in B(1)\stackrel{\triangle}{=}\{(x,y)^\top\in\mathbb R^2| x^2+y^2\leq 1\}$ a.e. $t\in(0,T)$. Set by $\phi_0(y_1(0),y_2(0),y_1(T),
y_2(T))=y_2(T)$ and $\psi(y_1(0),y_2(0),y_1(T),
y_2(T))=(y_1(0)-1,y_2(0))^\top$. Then, the optimal control problem is to minimize $\phi_0(y_1(0),y_2(0),y_1(T),
y_2(T))$, where $(y_1(\cdot),y_2(\cdot),u_1(\cdot),
u_2(\cdot))$ is subject to  (\ref{ex1}) and
$\psi(y_1(0),y_2(0),y_1(T),
y_2(T))=0$. Consider the control $\bar u(t)=(\bar u_1(t),\bar u_2(t))^\top\equiv (0,-1)^\top$. The corresponding trajectroy is
\begin{align*}
(\bar y_1(t),\bar y_2(t))=(1-t,-\frac{1}{3}t^3+3t^2-5t),\;\forall\,t\in[0,T].
\end{align*}
Then, we shall use \cite[Theorem 3.1 ]{gb}  and  Theorem \ref{c230}
 respectively to check whether $\bar u(\cdot)$ is optimal.

%

\end{exl}

{\it Solution.}\;It can be checked that $T_{B(1)}^\flat((0,-1)^\top)=\{(x,y)^\top| x\in\mathbb R, y\geq 0
\}$, which is strictly larger than $B(1)-(0,-1)^\top$. We also have $T_{B(1)}^{\flat(2)}((0,-1)^\top, (1,0)^\top)=\{(x,y)|x\in\mathbb R, y\geq \frac{1}{2}\}$. By  the \cite[Theorem 3.1(i) ]{gb}(or Theorem \ref{f}), there exists a unique $\ell=(\ell_0,\ell_1,\ell_2)^\top\in\mathbb R^3\setminus\{0\}$ (up to a positive factor) satisfying $\ell_0\leq 0$ and (\ref{ex2}) (or (\ref{c240})), and we also have $-\ell_2=\ell_0$ and $\ell_1=(6T-T^2)\ell_0$. Without loss of generality, we assume $\ell_0=-1$. We observe that, for all $v(\cdot)\in L^\infty(0,T;\mathbb R^2)\cap(\mathcal U-\{\bar u(\cdot)\})$, the relation $\int_0^T\nabla_uH[t,\ell](v(t))dt<0$ holds, which means that, \cite[Theorem 3.1(ii) ]{gb}(the second order necessary condition) can not be applied to $\bar u(\cdot)$.  However, by Theorem \ref{f} we know that $v(t)\equiv(1,0)^\top\in T_{B(1)}^\flat(\bar u(t))\setminus(B(1)-\{\bar u(t)\})$  $(t\in[0,T])$ is a singular direction. Set by $\sigma(t)\equiv (0,\frac{1}{2})^\top\in T_{B(1)}^{\flat(2)}(\bar u(t),(1,0)^\top)$ ($t\in[0,T]$).  The left hand side of (\ref{c244}) is reduced to
$T(-\frac{1}{3}T^2+\frac{5}{2}T+\theta-2)>0$, which contradicts (\ref{c244}). Consequently, $\bar u(\cdot)$ is not an optimal control.
$\Box$

\medskip

Then, we would apply Theorem \ref{c230} to the following problem
\begin{description}
\item[(OCPE)] Minimize $J(u(\cdot))=\int_0^Tf^0(t,y(t),u(t))dt$
over $u(\cdot)\in \mathcal U$ subject to (\ref{s2}), $y(0)=y_0$ and $y(T)=y_1$, 
where    $y_0, y_1\in M$ are fixed,  and $f:\mathbb R^+\times M\times U\to T M$ and $f^0:  \mathbb R^+\times M\times U\to  \mathbb R$ are given maps.
\end{description}
The Hamiltonian function $H^e:[0,T]\times T^*M\times U\times \mathbb R\to \mathbb R$ is given by
$
H^e(t,y,p,u,\ell)=p(f(t,y,u))+\ell f^0(t,y,u)
$, where $ (t,y,p,u,\ell)\in [0,T]\times T^*M\times U\times \mathbb R$.
We need the following assumption:
\begin{description}
\item[$(C_e)$] The maps $f$ and $f^0$ are measurable in $t$, $C^2$ in $(x,u)$, and there exist a constant $K>1$ and $x_0\in M$ such that 
\begin{align*}
&|f^0(s,x_1,u_1)-f^0(s,x_2,u_2)|\leq K(\rho(x_1,x_2)+|  u_1-u_2|),
\\ &|f^0(s,x_0,u_1)|\leq K,
\\& |\nabla_xf^0(s,x_1,u_1)-L_{\hat x_1 x_1}\nabla_xf^0(s,\hat x_1,u_2)|\leq K\big(\rho(x_1,x_2)+|u_1-u_2|\big),
\\&|\nabla_uf^0(s,x_1,u_1)-\nabla_uf^0(s,\hat x_1,u_2)|\leq K\big(\rho(x_1,x_2)+|u_1-u_2|\big),
\end{align*}
hold for all $x_1, \hat x_1,x_2\in M$, $u_1,u_2\in U$ and $s\in\mathbb R^+$, with $\rho(x_1,\hat x_1)<\min\{i(x_1),i(\hat x_1)\}$.
\end{description}

Then, the corresponding second order necessary condition is stated as follows.

\begin{cor}\label{c355} Assume that conditions $(C1)$ and $(C_e)$ hold,  that there exist a constant $K>1$ and $x_0\in M$ such that the first two lines of  (\ref{10}) and (\ref{116}) hold 
for all $x_1,\hat x_1\in M$ with $\rho(x_1,\hat x_1)<\min\{i(x_1),i(\hat x_1)\}$, $u_1,u_2\in U$ and $s\in\mathbb R^+$, and that
$(\bar u(\cdot), \bar y(\cdot))$ is an optimal pair of problem $(OCPE)$. Then, there exist $\ell_0\leq 0$ and $\varphi\in T_{y_1}^*M$ such that $(\ell_0,\varphi)\neq 0$ and
\begin{eqnarray}
\label{c340}
\nabla_uH^e[t,\ell_0,\varphi] (w(t))\leq 0,\quad \textrm{a.e.}\;t\in[0,T],
\end{eqnarray}
holds for any $w(\cdot)\in L^2(0,T;\mathbb R^m)$ with $w(t)\in T_U^\flat(\bar u(t))$ a.e. $t\in[0,T]$,
where we have used notion (\ref{434}), $p_{\ell_0\varphi}(\cdot)$ solves
\begin{equation}\label{c350}\left\{\begin{array}{l}
\nabla_{\dot{\bar y}(t)}p_{\ell_0\varphi}=-\nabla_xf[t](p_{\ell_0\varphi}(t),\cdot)-\ell_0\nabla_xf^0[t],\;\textrm{a.e.}\;t\in(0,T),
\\ p_{\ell_0\varphi}(T)=\varphi,
\end{array}\right.\end{equation} and we have adopted 
$[t,\ell_0,\varphi]=(t,\bar y(t),p_{\ell_0\varphi}(t),\bar u(t),\ell_0)$ for abbreviation. Moreover,   for any  $v(\cdot)\in   L^2(0,T;\mathbb R^m)$ with $v(t)\in T_U^{\flat}(\bar u(t))$ a.e. $t\in[0,T]$ and vector field  $X_v(\cdot)$ along $\bar y(\cdot)$ satisfying (\ref{c227}), $ X_v(0)=0$, $X_v(T)=0$,
and $
\int_0^T\left(\nabla_xf^0[t](X_v(t))+\nabla_uf^0[t](v(t))\right)dt\leq 0$, 
there exist $(\hat\ell_0,\hat\varphi)\in\big((-\infty,0]\times T_{y_1}^*M\big)\setminus\{0\}$  satisfying  (\ref{c340}) and (\ref{c350}) (with $(\ell_0,\varphi)$ replaced by $(\hat\ell_0,\hat\varphi)$), such that
\begin{equation}\label{c356}\begin{array}{l}
\int_0^T\nabla_u H^e[t,\hat\ell_0,\hat\varphi]\sigma(t)dt+\frac{1}{2}\int_0^T\left(\nabla_x^2H^e[t,\hat\ell_0,\hat\varphi](X_v(t),X_v(t))\right.
\\ +2\nabla_u\nabla_xH^e[t,\hat\ell_0,\hat\varphi](X_v(t),v(t))+\nabla_u^2H^e[t,\hat\ell_0,\hat\varphi](v(t),v(t))
\\\left.-R(\tilde p_{\hat\ell_0\hat\varphi}(t),X_v(t),f[t],X_v(t))\right)dt\leq 0,
\end{array}\end{equation}
holds for all $\sigma(\cdot)\in \mathcal B$,
where $p_{\hat \ell_0\hat\varphi}(\cdot)$ is the solution to (\ref{c350}) with $(\ell_0,\varphi)$ replaced by $(\hat\ell_0,\hat\varphi)$,  $\tilde p_{\hat\ell_0\hat\varphi}(t)$ ($t\in[0,T]$) is the dual vector of $p_{\hat\ell_0\hat\varphi}(t)$,  and for $(t,y,\eta,u,\ell)\in[0,T]\times T^*M\times U\times\mathbb R$,  the corresponding values of $\nabla_u {H^e}$, $\nabla_xH^e$,  $\nabla_x^2{H^e}$, $\nabla_u\nabla_x{H^e}$ and  $\nabla_u^2H^e$ at $(t,y,\eta,u,\ell)$   are respectively defined by
\begin{align*}
&\nabla_uH^e(t,y,\eta,u,\ell)(V)=\frac{d}{ds}\Big|_{s=0}H^e(t,y,\eta,u+sV,\ell),
\\
&\nabla_xH^e(t,y,\eta,u,\ell)(X)=\nabla_xf(t,y,u)(\eta,X)+\ell\nabla_xf^0(t,y,u)(X),
\\& \nabla_x^2H^e(t,y,\eta,u,\ell)(X,X)=\nabla_x^2f(t,y,u)(\eta,X,X)+\ell\nabla_x^2f^0(t,y,u)(X,X),
\\&\nabla_u\nabla_xH^e(t,y,\eta,u,\ell)(X,V)=\frac{d}{ds}\Big|_{s=0}\nabla_xH^e(t,y,\eta,u+sV,\ell)(X),
\\& \nabla_u^2H^e(t,y,\eta,u,\ell)(V,V)=\frac{d^2}{ds^2}\Big|_{s=0}H(t,y,\eta,u+sV,\ell),
\end{align*}
where $(X,V)\in T_yM\times\mathbb R^m$.
\end{cor}
By Theorem \ref{c230}, we can use the idea in \cite[Section 3.1]{dzgJ} to prove Corollary \ref{c355}, and we omit its proof.

For the case that $U$ is open (not necessarily convex), \cite[Theorem 3.2]{cdz} gives the second order necessary condition of optimal pairs of problem $(OCPE)$, which is in fact obtained through convex variation. Thus, Corollary  \ref{c355} is an complement to it.

\setcounter{equation}{0}
\hskip\parindent \section{Variations of Trajectories}
\def\theequation{3.\arabic{equation}}
 In this section,  we will compute variations of   (\ref{s2}), in the sense of convex variation.

\begin{pro}\label{439} Assume that  conditions $(C2)$ and $(C3)$ hold, and that $U\subset\mathbb R^m$ ($m\in\mathbb N$) is convex. Let $\bar u(\cdot)\in L^2(0,T;\mathbb R^m)\cap\mathcal U$ and $(\bar u(\cdot),\bar y(\cdot))$ satisfy (\ref{s2}).
Fix   $W \in T_{\bar y(0)}M$ and $v(\cdot)\in L^2(0,T;\mathbb R^m) $. Let $X_v(\cdot)$ be a vector field along $\bar y(\cdot)$ and satisfy (\ref{c227}).  Assume a set $\{\sigma_\epsilon(\cdot)\}_{\epsilon>0}\subset L^2(0,T;\mathbb R^m)$ is bounded in $L^2(0,T;\mathbb R^m)$, with $\sup_{\epsilon>0}\|\sigma_\epsilon(\cdot)\|_{L^2(0,T;\mathbb R^m)}=C_\sigma$.
   Denote by $y_\epsilon(\cdot)$ the solution to (\ref{s2}) corresponding to the control $u_\epsilon(\cdot ):=\bar u(\cdot)+\epsilon v(\cdot)+\epsilon^2\sigma_\epsilon(\cdot)$ and the initial state $y_\epsilon(0)=\exp_{\bar y(0)}(\epsilon X_v(0)+\epsilon^2W)$. Also denote  by   $Y_{\sigma_\epsilon W}^{X_v}(\cdot)$  the  solution to 
\begin{equation}\label{o2}
\left\{\begin{array}{l}\nabla_{\dot{\bar y}(t)}Y_{\sigma_\epsilon W}^{X_v}(Z)=\nabla_xf[t](Z,Y_{\sigma_\epsilon W}^{X_v}(t))+\nabla_uf[t](Z,\sigma_\epsilon(t))+\nabla_u\nabla_xf[t](Z,X_{v}(t),v(t))
\\\quad\quad\quad\quad\quad\quad-\frac{1}{2}R(\tilde Z,X
_{v}(t),\dot{\bar y}(t),X_{v}(t))+\frac{1}{2}\nabla_x^2f[t](Z,X_{v}(t),X_{v}(t))
\\ \quad\quad\quad\quad\quad\quad+\frac{1}{2}\nabla_u^2f[t](Z, v(t),v(t)),\quad a.e. \,t\in(0,T], \forall\; Z\in T^*M,\cr Y_{\sigma_\epsilon W}^{X_v}(0)=W,
\end{array}\right.\end{equation}
where we adopt notion (\ref{434}), and 
\begin{align*}
&\nabla_u^2f[t](Z,v(t),v(t))=\frac{\partial^2}{\partial s^2}\big|_{s=0}f\big(t,\bar y(t),\bar u(t)+sv(t)\big)(Z),
\\&\nabla_u\nabla_xf[t]\big(Z,X_{v}(t),v(t)\big)=\frac{\partial}{\partial s}\big|_{s=0}\nabla_xf\big(t,\bar y(t),\bar u(t)+sv(t)\big)\big(Z,X_{v}(t)\big).
\end{align*}
 Then, for any $\alpha>0$, there exists $\epsilon^0>0$ such that
\begin{equation}\label{427}
|V_\epsilon(t)-\epsilon X_{v}(t)-\epsilon^2Y_{\sigma_\epsilon W}^{X_v}(t)|\leq \alpha \epsilon^2,\quad\forall\; t\in[0,T],\,\forall\,\epsilon\in[0,\epsilon^0],\end{equation}
where
\begin{equation}\label{438}
V_\epsilon(t):=\exp_{\bar y(t)}^{-1}y_\epsilon(t),\quad t\in[0,T].
\end{equation}
\end{pro}

\textbf{Proof.}
The proof is split into three steps.

{\it Step 1.}  We claim that, there exists $\hat\epsilon>0$ such that (\ref{438}) can be defined for $\epsilon\in[0,\hat\epsilon]$, and 
\begin{align}\label{ccs15}\begin{array}{lll}
|V_\epsilon(t)|&=&\rho(y_\epsilon(t),\bar y(t))
\\&\leq &\epsilon e^{\frac{1}{2}(K+\frac{1}{2})T}\Big(|X_v(0)+\epsilon W|^2+4K^2\|v\|^2_{L^2(0,T;\mathbb R^m)}+4K^2C^2_\sigma\Big)^{\frac{1}{2}},\,t\in[0,T],
\end{array}
\end{align}
for all $\epsilon\in[0,\hat\epsilon]$.

In fact, let $\epsilon_0\in(0,1]$ be such that $\epsilon_0|X_v(0)|+\epsilon_0^2|W|<i(\bar y(0))$, where the injectivity radius $i(y_0)$ of the point $y_0$ is defined in \cite[Section 2.1]{cdz}.  Then, by \cite[Lemma 5.2]{cdz}, the triangle inequality of $\rho(\cdot,\cdot)$ and \cite[Lemma 2.2]{cdz}, we have 
\begin{align}\label{ccs13}\begin{array}{ll}
&\rho(y_\epsilon(t),\bar y(0))\leq \rho(y_\epsilon(t),y_\epsilon(0))+\rho(y_\epsilon(0),\bar y(0))
\\\leq & (1+\rho(x_0,\bar y(0))+\epsilon |X_v(0)|+\epsilon^2|W|)e^{Kt},\quad\forall\,t\in[0,+\infty),\quad\forall\,\epsilon\in[0,\epsilon_0].
\end{array}\end{align}
By HopfRinow theorem (see \cite[Theorem. 16, p. 137]{p1}), there exists $\delta>0$ such that $i(y)\geq \delta$ for all $y\in M$ with $\rho(\bar y(0),y)\leq (1+\rho(x_0,\bar y(0))+\epsilon_0 |X_v(0)|+\epsilon_0^2|W|)e^{KT}$.  	Similar to (\ref{ccs13}), for any $t, \hat t\in[0,+\infty)$, we have
\begin{align}\label{ccs19}
\rho(y_\epsilon(t),\bar y(t))\leq \Big(2+2\rho(x_0,\bar y(0))+\epsilon|X_v(0)|+\epsilon^2|W|\Big)(e^{Kt}-e^{K\hat t})+\rho(y_\epsilon(\hat t),\bar y(\hat t)),
\end{align}
for all $\epsilon\in [0,\epsilon_0]$.

Let $\hat t=0$ in (\ref{ccs19}).  We take $\epsilon_1>0$ such that $\epsilon_1|X_v(0)|+\epsilon_1^2|W|=\frac{\delta}{2}$, and take $T_1>0$ such that $(2+2\rho(x_0,\bar y(0))+\epsilon_1|X_v(0)|+\epsilon_1^2|W|\Big)(e^{KT_1}-1)=\frac{\delta}{2}$. Thus, we can define (\ref{438}) over $[0,T_1]\cap [0,T]$, and obtain   by \cite[(2.17) \& (2.19)]{cdz}, \cite[(2.2)]{dzgJ} and $(C2)$ that
\begin{align}\label{esd}\begin{array}{lll}
\frac{d}{dt}\rho^2(y_\epsilon(t),\bar y(t))&=&\nabla_1\rho^2(y_\epsilon(t),\bar y(t))\Big(f(t,y_\epsilon(t),u_\epsilon(t))-L_{\bar y(t)y_\epsilon(t)}f(t,\bar y(t),u_\epsilon(t))\Big)
\\[2mm]&&+\nabla_2\rho^2(y_\epsilon(t),\bar y(t))\Big(f(t,\bar y(t),u_\epsilon(t))-f[t]\Big)
\\[2mm]&\leq & (K+\frac{1}{2})\rho^2(y_\epsilon(t),\bar y(t))+4\epsilon^2K^2(|v(t)|^2+|\sigma_\epsilon(t)|^2),
\end{array}
\end{align}
for all $t\in[0,T_1]\cap [0,T]$ and $\epsilon\in[0,\epsilon_1]$. 
 By the Gronwall's inequality and \cite[Lemma 2.2]{cdz}, we obtain that  (\ref{ccs15}) holds
for all $t\in[0, T_1]\cap [0,T]$ and $\epsilon\in[0,\epsilon_1]$.

 If $T_1<T$,  we set
 $\hat t=T_1$ in (\ref{ccs19}),   $\hat\epsilon=\min\{\epsilon_1, \frac{\delta}{2 C_{v,\sigma}}\}$ with $(C_{v,\sigma})^2=e^{(K+\frac{1}{2})T}\Big(2|X_v(0)|^2+2\epsilon_1^2|W|^2+4K^2\|v\|^2_{L^2(0,T;\mathbb R^m)}+4K^2C^2_\sigma\Big)$, and $T_2>T_1$ satisfying $\Big(2+2\rho(x_0,\bar y(0))+\hat\epsilon|X_v(0)|+\hat\epsilon^2|W|\Big)(e^{KT_2}-e^{KT_1})=\frac{\delta}{2}$. Then we can define (\ref{438}) over $[0,T_2]\cap [0,T]$ for $\epsilon\in[0,\hat\epsilon]$, and consequently (\ref{esd}) holds for $t\in[0,T_2]\cap[0,T]$ and $\epsilon\in[0,\hat\epsilon]$. Analogously we get (\ref{ccs15}) for $t\in[0,T_2]\cap [0,T]$ and $\epsilon\in[0,\hat\epsilon]$.
Recursively, if $T_i<T$, we take $\hat t=T_i$ ($i\geq 2$), and set $T_{i+1}>T_i$  such that 
\begin{align}\label{ccs30}
\Big(2+2\rho(x_0,\bar y(0))+\hat\epsilon|X_v(0)|+\hat\epsilon^2|W|\Big)(e^{KT_{i+1}}-e^{KT_i})=\frac{\delta}{2}.
\end{align}
 Then, one can define (\ref{438}) over $[0,T_{i+1}]\cap [0,T]$. Consequently, we obtain (\ref{ccs15}) for $t\in[0,T_{i+1}]\cap [0,T]$ and $\epsilon\in[0,\hat\epsilon]$.  It follows from (\ref{ccs30}) that there eixsts $I\in\mathbb N$, such that $T_I>T$, and consequently,  (\ref{ccs15}) holds for $t\in[0,T]$ and $\epsilon\in[0,\hat\epsilon].$

{\it Step 2.} Let $\{e_1,\cdots, e_n\}\subset T_{\bar y(0)}M$ be an orthonormal basis at $\bar y(0)$, i.e. $\langle e_i,e_j\rangle=\delta_i^j$ for $i,j=1,\cdots,n$, where $\delta_i^j$ is the usual Kronecker symbol.  Denote by $\{d_i\}_{i=1}^n\subset  T_{\bar y(0)}^*M$ the dual basis to it. 
For $t\in[0,T]$, set respectively by $e_i(t)=L^{\bar y(\cdot)}_{\bar y(0)\bar y(t)}e_i$ and $d_i(t)=L^{\bar y(\cdot)}_{\bar y(0)\bar y(t)}d_i$ for $i=1,\cdots,n$, where $L^{\bar y(\cdot)}_{\bar y(0)\bar y(t)} $ is the parallel translation along the curve $\bar y(\cdot)$ and from $\bar y(0)$ to $\bar y(t)$, see \cite[Section 2.2]{cdz} for its detailed definition, and then $\nabla_{\dot{\bar y}(t)}e_i(\cdot)=0$ and $\nabla_{\dot{\bar y}(t)}d_i(\cdot)=0$ for $i=1,\cdots,n$.
We deduce from \cite[(2.7)\&(2.6)]{cdz} that, $\{e_i(t)\}_{i=1}^n\subset T_{\bar y(t)}M$ is an orthonormal basis, and  $\{d_i(t)\}_{i=1}^n\subset T_{\bar y(t)}^*M$ is the dual basis to it.

For $\epsilon\in[0,\hat\epsilon]$ and $t\in[0,T]$, it follows from "Step 1", \cite[Lemma 2.1]{cdz} and the definition of exponential map (see \cite[Section 2.1]{cdz})  that,  there exists a unique geodesic 
\begin{align}\label{ccs100}
\beta(\theta;t)=\exp_{\bar y(t)}(\theta V_\epsilon(t)),\quad \forall \,\theta\in[0,1],
\end{align}
connecting $\beta(0;t)=\bar y(t)$ and $\beta(1;t)=y_\epsilon(t)$.
For $\theta\in[0,1]$, denote by $L_{\bar y(t)\beta(\theta;t)}: T_{\bar y(t)}M\to T_{\beta(\theta;t)}M$ the parallel translation along the geodesic $\beta(\cdot;t)$ ($t\in[0,T]$ is fixed), from $\bar y(t)$ to $\beta(\theta;t)$.  By \cite[(2.6)]{cdz} we know that $\{e_i(t)\}_{i=1}^n$ and $\{L_{\bar y(t)\beta(\theta;t)} e_i(t)\}_{i=1}^n$ are respectively  orthonormal bases at $T_{\bar y(t)}M$ and $T_{\beta(\theta;t)}M$ for each $t\in[0,T]$ and $\theta\in[0,1]$.
Thus, we can write
\begin{align}\label{ccs35}
V_\epsilon(t)=\sum_{i=1}^na_i^\epsilon(t)e_i(t),\;t\in[0,T],
\end{align}
where $\{a_1^\epsilon(t),\cdots, a_n^\epsilon(t)\}$ satisfies
$
\sum_{i=1}^na_i^\epsilon(t)^2=|V_\epsilon(t)|^2,\;\forall\,t\in[0,T].
$
By the linearity of $L_{\bar y(t)\beta(\theta;t)}: T_{\bar y(t)}M\to T_{\beta(\theta;t)}M$ and \cite[(2.8)]{cdz} we have
\begin{align}\label{ccs36}
L_{\bar y(t)\beta(\theta;t)}V_\epsilon(t)=\sum_{i=1}^na_i^\epsilon(t)L_{\bar y(t)\beta(\theta;t)} e_i(t),\;|L_{\bar y(t)\beta(\theta;t)}V_\epsilon(t)|=|V_\epsilon(t)|,\;\forall\,t\in[0,T].
\end{align}

Fix any $t\in[0,T]$ and $i\in\{1,\cdots,n\}$.    By \cite[Lemma 2.2]{cdz}, Newton-Leibniz formula and the exchange of integral order, we derive
\begin{equation}\label{l30}\begin{array}{ll}
&\langle\nabla_{\dot{\bar y}(t)}V_\epsilon,e_i(t)\rangle
\\=& \frac{d}{dt}\langle V_\epsilon,e_i(t)\rangle
\\ =& -\frac{d}{dt}\Big(\frac{1}{2}\nabla_1\rho^2(\bar y(t),y_\epsilon(t))(d_i(t))\Big)
\\:=& P_1^i(t)+\nabla_uf[t]\Big(d_i(t),\epsilon v(t)+\epsilon^2\sigma_\epsilon(t)\Big)
\\&+\int_0^1\nabla_u^2 f\{t,\tau\}_\epsilon\Big(d_i(t),\epsilon v(t)+\epsilon^2\sigma_\epsilon(t), \epsilon v(t)+\epsilon^2\sigma_\epsilon(t)\Big)(1-\tau)d\tau,
\end{array}\end{equation}
where
\begin{equation}\label{o3}
\{t,\tau\}_\epsilon:=(t,\bar y(t),\bar u(t)+\tau\epsilon v(t)+\tau\epsilon^2\sigma_\epsilon(t)),\quad\forall\, t\in[0,T],\; \tau\in[0,1],
\end{equation}
and
$$\begin{array}{ll}
&P_1^i(t)
\\:=&-\frac{1}{2}\Big\{\nabla_2\nabla_1\rho^2(\bar y(t),y_\epsilon(t))(e_i(t),f(t,y_\epsilon(t),u_\epsilon(t)))
\\&-\nabla_2\nabla_1\rho^2(\bar y(t),\bar y(t))(e_i(t),f(t,\bar y(t),u_\epsilon(t)))
\\&+\nabla_1^2\rho^2(\bar y(t),y_\epsilon(t))(e_i(t),f[t])-\nabla_1^2\rho^2(\bar y(t),\bar y(t))(d_i(t),f[t])\Big\}.
\end{array}$$

From the definitions of geodesic (see \cite[(2.1)]{cdz}) and parallel translation (see \cite[Section 2.2]{cdz}) and \cite[Lemma 2.1]{cdz}, we obtain
\begin{align}\label{ccs34}
\frac{\partial}{\partial \theta}\beta(\theta;t)=L_{\bar y(t)\beta(\theta;t)}\frac{\partial}{\partial \theta}\Big|_0\beta(\theta;t)=L_{\bar y(t)\beta(\theta;t)} V_\epsilon(t),\quad\theta\in[0,1].
\end{align}
Applying \cite[Lemma 2.2]{cdz} and Newton-Leibniz formula to $P_1^i(t)$, we have
\begin{align*}
P_1^i(t)=I^i_1(t)+I^i_2(t)+I^i_3(t),
\end{align*}
where
\begin{align*}
I_1^i(t)=&-\frac{1}{2}\int_0^1\Big[\nabla_2\nabla_1^2\rho^2
(\bar y(t),\beta(\theta;t))\Big(e_i(t),f[t],\frac{\partial}{\partial\theta}\beta(\theta;t)\Big)
\\&-\nabla_2\nabla_1^2\rho^2
(\bar y(t),\beta(0;t))\Big(e_i(t),f[t],\frac{\partial}{\partial\theta}\Big|_0\beta(\theta;t)\Big)\Big]d\theta,
\end{align*}
\begin{align*}
I_2^i(t)=&-\frac{1}{2}\int_0^1\Big[\nabla_2^2\nabla_1
\rho^2(\bar y(t),\beta(\theta;t))\Big(e_i(t),f(t,\beta(\theta;t),u_\epsilon(t)),\frac{\partial}{\partial\theta}\beta(\theta;t)\Big)
\\&-\nabla_2^2\nabla_1\rho^2(\bar y(t),\beta(0;t))\Big(e_i(t),f(t,\beta(0;t),u_\epsilon(t)),\frac{\partial}{\partial\theta}\Big|_0\beta(\theta;t)\Big)\Big]d\theta,
\end{align*}
and
\begin{align*}
I_3^i(t)=&\nabla_xf[t](d_i(t),V_\epsilon(t))+\int_0^1\nabla_u\nabla_xf(t,\bar y(t), u_\epsilon^\theta(t))\Big(d_i(t),V_\epsilon(t), \epsilon v(t) 
\\&+\epsilon^2\sigma_\epsilon(t)\Big)d\theta
-\frac{1}{2}\int_0^1\Big[\nabla_2\nabla_1\rho^2(\bar y(t),\beta(\theta;t))\Big(e_i(t),\nabla_{\frac{\partial}{\partial\theta}\beta(\theta;t)}f(t,\cdot,u_\epsilon
(t))\Big)
\\&-\nabla_2\nabla_1\rho^2(\bar y(t),\bar y(t))\Big(e_i(t),\nabla_{\frac{\partial}{\partial\theta}|_0\beta(\theta;t)}f(t,\cdot,u_\epsilon(t))\Big)\Big]d\theta,
\end{align*}
with $u_\epsilon^\theta(t)=\bar u(t)+\theta(\epsilon v(t)+\epsilon^2\sigma_\epsilon(t)), \theta\in[0,1]$. We use Newton-Leibniz formula
again to the above three items, exchange the integration order, and get
\begin{align}\label{ccs50}
&I_1^i(t)
=-\frac{1}{4}\nabla_2^2\nabla_1^2\rho^2(\bar y(t),\bar y(t))\Big(e_i(t),f[t],V_\epsilon(t),V_\epsilon(t)\Big)+\hat I_1^i(t),
\\\label{ccs51}
&I_2^i(t)
=-\frac{1}{4}\nabla_2^3\nabla_1\rho^2(\bar y(t),\bar y(t))\Big(e_i(t),f[t],V_\epsilon(t),V_\epsilon(t)\Big)+\hat I_2^i(t),
\end{align}
and 
\begin{align}\label{ccs52}
\begin{array}{ll}&I_3^i(t)
\\=&\nabla_xf[t](e_i(t),V_\epsilon(t))+\epsilon\nabla_u\nabla_xf[t](d_i(t),V_\epsilon(t),v(t))
\\&-\frac{1}{4}\nabla_2\nabla_1\rho^2(\bar y(t),\bar y(t))\Big(e_i(t),\nabla_{\frac{\partial}{\partial\tau}|_0\beta(\tau;t)}\nabla_{\frac{\partial}{\partial\tau}\beta(\tau;t)}f(t,\cdot,\bar u(t))\Big)+\hat I_3^i(t),
\end{array}\end{align}
where we have used (\ref{ccs35}), (\ref{ccs36}), (\ref{ccs34}), and the fact that $\beta(\cdot;t)$ is a geodesic, 
\begin{align*}
\hat I_1^i(t)=&-\frac{1}{2}\sum_{k,l=1}^n\int_0^1\Big[\nabla_2^2\nabla_1^2
\rho^2(\bar y(t),\beta(\theta;t))\Big(e_i(t),f[t],L_{\bar y(t)\beta(\theta;t)}e_k(t),L_{\bar y(t)\beta(\theta;t)}e_l(t)\Big)
\\&-\nabla_2^2\nabla_1^2
\rho^2(\bar y(t),\bar y(t))\Big(e_i(t),f[t],e_k(t),e_l(t)\Big)\Big](1-\theta)d\theta\, a_k^\epsilon(t)a_l^\epsilon(t),
\end{align*}
\begin{align*}
\hat I_2^i(t)=&-\frac{1}{4}\nabla_2^3\nabla_1\rho^2(\bar y(t),\bar y(t))\Big(e_i(t),
f(t,\bar y(t),u_\epsilon(t))-f[t],V_\epsilon(t),V_\epsilon(t)\Big)
\\&-\sum_{k,l=1}^n\frac{1}{2}\int_0^1\Big[\nabla_2^3\nabla_1
\rho^2
(\bar y(t),\beta(\tau;t))\Big(e_i(t),
f(t,\beta(\tau;t),u_\epsilon(t)),
\\&L_{\bar y(t)\beta(\theta;t)}e_k(t),L_{\bar y(t)\beta(\theta;t)}e_l(t)\Big)-\nabla_2^3\nabla_1\rho^2(\bar y(t),\bar y(t))\Big(e_i(t),
f(t,\bar y(t),u_\epsilon(t)),
\\&e_k(t),e_l(t)\Big)\Big]a_k^\epsilon(t)a_l^\epsilon(t)(1-\tau)d\tau
-\frac{1}{
2}\int_0^1\nabla_2^2\nabla_1\rho^2
(\bar y(t),\beta(\tau;t))\Big(e_i(t),\\&\nabla_{\frac{\partial}{\partial\tau}\beta(\tau;t)}f(t,\cdot,u_\epsilon(t)),\frac{\partial}{\partial\tau}\beta(\tau;t)\Big)(1-\tau)d\tau,
\end{align*}
and
\begin{align*}
&\hat I_3^i(t)
\\=
&-\frac{1}{4}\nabla_2\nabla_1\rho^2(\bar y(t),\bar y(t))\Big(e_i(t),\nabla_{\frac{\partial}{\partial\tau}|_0\beta(\tau;t)}\nabla_{\frac{\partial}{\partial\tau}\beta(\tau;t)}[f(t,\cdot,\bar u(t))-f(t,\cdot, u_\epsilon(t))]\Big)
\\&-\frac{1}{2}\int_0^1\Big[\nabla_2\nabla_1\rho^2(\bar y(t),\beta(\tau;t))\Big(e_i(t),\nabla_{\frac{\partial}{\partial\tau}\beta(\tau;t)}\nabla_{\frac{\partial}{\partial\tau}\beta(\tau;t)}f(t,\cdot,u_\epsilon(t))\Big)
\\&-\nabla_2\nabla_1\rho^2(\bar y(t),\bar y(t))\Big(e_i(t),\nabla_{\frac{\partial}{\partial\tau}|_0\beta(\tau;t)}\nabla_{\frac{\partial}{\partial\tau}\beta(\tau;t)}f(t,\cdot,u_\epsilon(t))\Big)\Big](1-\tau)d\tau
\\&+\epsilon\int_0^1\Big[\Big(\nabla_u\nabla_xf(t,\bar y(t),u_\epsilon^\theta(t))-\nabla_u\nabla_xf[t]\Big)\Big(d_i(t),V_\epsilon(t),v(t)\Big)
\\&+\epsilon\nabla_u\nabla_xf(t,\bar y(t),u_\epsilon^\theta(t))\Big(d_i(t),V_\epsilon(t),\sigma_\epsilon(t)\Big)\Big]d\theta
\\&-\frac{1}{2}\int_0^1\nabla_2^2\nabla_1\rho^2(\bar y(t),\beta(\tau;t))\Big(e_i(t),\nabla_{\frac{\partial}{\partial\tau}\beta(\tau;t)}f(t,\cdot,u_\epsilon(t)),
\\&\frac{\partial}{\partial\tau}\beta(\tau;t)\Big)(1-\tau)d\tau.
\end{align*}
Following the same argument as that used in \cite[(5.18)]{cdz}, we have
\begin{align}\label{ccs53}\begin{array}{ll}
&\langle Z(\beta(s;t)),\nabla_{\frac{\partial}{\partial\tau}|_s\beta(\tau;t)}\nabla_{\frac{\partial}{\partial\tau}\beta(\tau;t)}f(t,\cdot,u)\rangle
\\[2mm]=&\nabla^2 f(t,\beta(s;t),u)\Big(\tilde Z,\frac{\partial}{\partial\tau}\Big|_s\beta(\tau;t),\frac{\partial}{\partial\tau}\Big|_s\beta(\tau;t)\Big),
\end{array}
\end{align}
where $(s,u,Z)\in[0,1]\times U\times T M$, and $\tilde Z$ is the dual covector of $Z$.

{\it Step 3.} 
Recall (\ref{c227}), (\ref{o2}), (\ref{l30}), (\ref{ccs50}) - (\ref{ccs52}). For $\epsilon\in[0,\hat\epsilon]$, by (\ref{ccs53}), \cite[Lemma 2.2 \& Lemma 2.3]{cdz} and Newton-Leibniz formula, we get
\begin{align}\label{ccs41}\begin{array}{ll}
&\langle V_\epsilon(t)-\epsilon X_{v}(t)-\epsilon^2Y_{\sigma_\epsilon W}^{X_v}(t),e_i(t)\rangle
\\[2mm]=&\displaystyle\int_0^t\langle\nabla_{\dot{\bar y}(s)}V_\epsilon-\epsilon \nabla_{\dot{\bar y}(s)}X_{v}-\epsilon^2\nabla_{\dot{\bar y}(s)}Y_{\sigma_\epsilon W}^{X_v},e_i(s)\rangle ds
\\[2mm]=&\displaystyle\int_0^t\Big\{\nabla_xf[s]\Big(d_i(s),V_\epsilon(s)-\epsilon X_{v}(s)-\epsilon^2Y_{\sigma_\epsilon W}^{X_v}(s)\Big)
\\[2mm]&+\epsilon^2\displaystyle\int_0^1\Big[\nabla_u^2f\{s,\tau\}_\epsilon\Big(d_i(s),v(s),v(s)\Big)-\nabla_u^2f[s]\Big(d_i(s),v(s),v(s)\Big)\Big](1-\tau)d\tau
\\[2mm]&-\frac{1}{2}R(e_i(s),V_\epsilon(s),f[s],V_\epsilon(s))+\frac{\epsilon^2}{2}R(e_i(s),X_{v}(s),f[s],X_{v}(s))
\\[2mm]&+\frac{1}{2}\nabla_x^2f[s]\Big(d_i(s),V_\epsilon(s),V_\epsilon(s)\Big)-\frac{\epsilon^2}{2}\nabla_x^2f[s]\Big(d_i(s),X_{v}(s),X_{v}(s)\Big)
\\[2mm]&+\epsilon\nabla_u\nabla_xf[s]\Big(d_i(s),V_\epsilon(s)-\epsilon X_{v}(s),v(s)\Big)\Big\}ds+\epsilon^2P_2^i(t),
\end{array}\end{align}
where 
\begin{align*}
&P_2^i(t)
\\=&\int_0^t\Big\{\int_0^1\Big[
2\epsilon\nabla_u^2f\{s,\tau\}_\epsilon\Big(d_i(s),v(s),\sigma_\epsilon(s)\Big)
\\&+\epsilon^2\nabla_u^2f\{s,\tau\}_\epsilon\Big(d_i(s),\sigma_\epsilon(s),\sigma_\epsilon(s)\Big)\Big]
(1-\tau)d\tau 
+\frac{1}{\epsilon^2}\Big(\hat I_1^i(s)+\hat I_2^i(s)+\hat I_3^i(s)\Big)\Big\}ds.
\end{align*}

It follows from $(C2)-(C4)$ and (\ref{ccs15}) that $|P_1^i(t)|$ is bounded for $t\in[0,T]$. By applying Gronwall's  inequality to (\ref{ccs41}), we obtain that, there exists a positive constant $C$ such that
\begin{align*}
|V_\epsilon(t)-\epsilon X_{v}(t)|\leq C\epsilon^2,\;\forall\,t\in[0,T],\,\forall\,\epsilon\in[0,\hat\epsilon].
\end{align*}
Applying Gronwall's inequality again to (\ref{ccs41}), we obtain from the above inequality, Lebesgue’s dominated convergence theorem, (\ref{ccs15}) and $(C2)-(C4)$ that, given any $\alpha>0$, there exists $\epsilon^0\in(0,\hat\epsilon]$ such that (\ref{427}) holds. The proof is concluded.
$\Box$

\setcounter{equation}{0}

\section{Proof of Theorem \ref{c230}}
\def\theequation{4.\arabic{equation}}
In Section  \ref{soop}, we obtain the second order necessary condition of an optimization problem (problem $(OP)$), see Theorem \ref{c126}. In Section \ref{psoc}, we transform problem $(OCP)$ into an optimization problem, which is a special case of problem $(OP)$, and prove Theorem \ref{c230} by Theorem \ref{c126}.
\subsection{An optimization problem}\label{soop}
Let $\mathcal X$ be a Banach space,  and $E\subset \mathcal X$ be a convex subset of it. Given  maps $\hat\phi_i:\mathcal X\to \mathbb R$ with $i=0,\cdots,j$, and $\hat\psi=(\hat\psi_1,\cdots,\hat\psi_k)^\top:\mathcal X\to \mathbb R^k$ ($k\in\mathbb N\cup\{0\}$), consider the following optimization problem. 
\begin{description}
\item[(OP)]  Find $\bar e\in E$ such that it minimizes $\hat\phi_0(e)$ with $e\in E$ subject to
\begin{equation}\label{c5}\begin{array}{l}
 \hat\phi_i(e)\leq 0,\;i=1,\cdots,j,
\\\hat\psi(e)=0.
\end{array}
\end{equation}
$\bar e$ is called a solution or a minimizer of problem $(OP)$.
\end{description}

Set $\hat\Phi:\mathcal X\to \mathbb R^{1+j+k}$ by
$\hat\Phi=(\hat\phi_0,\hat\phi_1,\cdots,\hat\phi_j,\hat\psi_1,\cdots,\hat\psi_k)^\top$.
Given any index set $I\subseteq\{0,1,\cdots,j\}$, we denote by $\hat\Phi_I=(\bar\phi_0,\bar \phi_1,\cdots,\bar\phi_j,\hat\psi_1,\cdots,\hat\psi_k)^\top$, where $\bar\phi_i=\hat\phi_i$ if $i\in I$, and $\bar\phi_i=0$ if $i\notin I$.

Assume $\bar e\in E$ is a minimizer of problem $(OP)$. 
Set by 
$$\begin{array}{l}
\hat I_{A}\equiv\{i\in\{1,\cdots,j\};\hat\phi_i(\bar e)=0\}\cup\{0\},
\\  \hat I_N\equiv\{i\in\{1,\cdots,j\}; \hat\phi_i(\bar e)<0\}.
\end{array}$$
To state the necessary condition of problem $(OP)$, we introduce the following condition.

\begin{description}
\item[(C5)]  $\hat\Phi$ is  Fr$\acute{\textrm{e}}$chet differentiable at $\bar e\in E$, and we denote its Fr$\acute{\textrm{e}}$chet  derivative
at $\bar e\in E$  by $D\hat\Phi(\bar e)=\big(D\hat\phi_0(\bar e),\cdots,D\hat\phi_j(\bar e), D\hat\psi_1(\bar e),\cdots,D\hat\psi_k(\bar e)\big)^\top$, where $D\hat\phi_i(\bar e)$ ($i=0,1,\cdots,j$) and $D\hat\psi_l(\bar e)$ ($l=1,\cdots,k$) are the Fr$\acute{\textrm{e}}$chet derivatives of $\hat\phi_i$ and $\hat\psi_l$ at $\bar e$ respectively. For any $x\in\mathcal X$, $$D\hat\Phi(\bar e)(x)\equiv\big(D\hat\phi_0(\bar e)(x),\cdots,D\hat\phi_j(\bar e)(x),D\hat\psi_1(\bar e)(x),\cdots,D\hat\psi_k(\bar e)(x)\big)^\top.$$   For each $y\in \mathcal X$, there exists 
$$D^2\hat\Phi(\bar e)(y)=\big(D^2\hat\phi_0(\bar e)(y),\cdots,D^2\hat\phi_j(\bar e)(y), D^2\hat\psi_1(\bar e)(y),\cdots,D^2\hat\psi_k(\bar e)(y)\big)^\top\in \mathbb R^{1+j+k},$$
  such that the following relation holds: for any  $\alpha>0$ and $C>0$, there exists $\epsilon_0>0$ depending on $\alpha$ and $C$, such that
\begin{align*}
|\hat\Phi(\bar e+\epsilon y+\epsilon^2\eta)-\hat\Phi(\bar e)-\epsilon D\hat\Phi(\bar e)(y)-\epsilon^2D\hat\Phi(\bar e)(\eta)-\frac{1}{2}\epsilon^2 D^2\hat\Phi(\bar e)(y)|\leq\alpha \epsilon^2,
\end{align*}
for all $\eta\in\mathcal X$ with $|\eta|\leq C$ and $\epsilon\in[0,\epsilon_0]$.

\end{description}

\begin{lem}\label{c115}  Assume that $\bar e\in E$ is a minimizer of $(OP)$ with $\hat\phi_0(\bar e)=0$, and that  $(C5)$ holds.  Let    $y\in T_E^\flat(\bar e)$ satisfy
\begin{equation}\label{c117}
\left\{\begin{array}{l}D\hat\phi_i(\bar e)( y)\leq 0,\; \textrm{for}\; i\in \hat I_{A},\\ D\hat\psi(\bar e)( y)=0,\quad T_E^{\flat(2)}(\bar e,y)\not=\emptyset,\end{array}\right.
\end{equation}
where $D\hat\psi(\bar e)=\big(D\hat\psi_1(\bar e),\cdots, D\hat\psi_k(\bar e)\big)^\top$ is a linear map from $\mathcal X$ to $\mathbb R^{k}$.
Set
$$\begin{array}{l}\hat I_0^\prime\equiv\hat I_N\cup \{i\in \hat I_{A}; D\hat\phi_i(\bar e)(y)<0\},
\\ \hat I_0^{\prime\prime}\equiv \{0,1,\cdots,j\}\setminus\hat I_0^\prime.
\end{array}$$
 Denote by 
\begin{equation}\label{c116}
\mathcal K\equiv\{D\hat\Phi_{\hat I_0^{\prime\prime}}(\bar e)(x)+\frac{1}{2}D^2\hat\Phi_{\hat I_0^{\prime\prime}}(\bar e)(y);x\in T_E^{\flat(2)}(\bar e, y)\}\subset \mathbb R^{1+j+k}
\end{equation}
and 
\begin{align*}
\mathcal K^{\hat\psi}=\{D\hat\psi(\bar e)(x)+\frac{1}{2}D^2\hat\psi(\bar e)(y); x\in T_E^{\flat(2)}(\bar e,y)\}\subset\mathbb R^k.
\end{align*}
Then $\mathcal K$ and $\mathcal K^{\hat\psi}$ are both convex. 
Moreover, we set by 
$
Y\equiv (Y_0,Y_1,\cdots,Y_j)^\top,
$ with
$$
Y_i=\left\{\begin{array}{l}D\hat\phi_i(\bar e)(y),\;i\in \hat I_{A},\cr
0,\quad \quad \quad\;i\notin \hat I_{A},\end{array}\right.
$$
 and by
\begin{equation}\label{c201}
Z\equiv (-\infty,0)^{j+1}-\{\lambda(\hat\phi(\bar e)+Y);\,\lambda>0\},
\end{equation}
where $\hat\phi=(\hat\phi_0,\hat\phi_1,\cdots,\hat\phi_j)^\top$. If $\mathcal K$ and $Z\times \{0\}$ can not be seperated by any linear functional: there does not exist
 $\ell\in\mathbb R^{1+j+k}\setminus\{0\}$ such that  
\begin{equation}\label{c129}
\ell^\top (D\hat\Phi_{\hat I_0^{\prime\prime}}(\bar e)(x)+\frac{1}{2}D^2\hat\Phi_{\hat I_0^{\prime\prime}}(\bar e)(y))\leq \ell^\top(z^\top,0)^\top,\;\forall x\in T_E^{\flat(2)}(\bar e,y),\;\forall z\in Z,
\end{equation} then, $\textrm{aff}\,\mathcal K^{\hat\psi}$ is a subspace of $\mathbb R^k$, where $\textrm{aff}\,\mathcal K^{\hat\psi}$ is the affine hull of $\mathcal K^{\hat\psi}$ (see \cite[p.6]{rock}), and the dimension  of $\textrm{aff}\,\mathcal K^{\hat\psi}$
 is bigger than  zero (see \cite[p. 4]{rock}). We denote it by $D(\mathcal K^{\hat\psi})$ . Moreover, 
 there exist $h_1,\cdots,h_{D(\mathcal K^{\hat\psi})+1}\in T_E^{\flat(2)}(\bar e,y)$ and $\delta_0>0$, such that
\begin{eqnarray}\label{c119}
B_{\textrm{aff}\,\mathcal K^{\hat\psi}}(\delta_0)\subseteq Int\, co\{D\hat\psi(\bar e)(h_l)+\frac{1}{2}D^2\hat\psi(\bar e)(y)\}_{l=1}^{D(\mathcal K^{\hat\psi})+1},
\\\label{c200}D\hat\phi_i(\bar e)(h_l)+\frac{1}{2}D^2\hat\phi_i(\bar e)(y)<0,\;l=1,\cdots,D(\mathcal K^{\hat\psi})+1,\;\textrm{if}\;i\in \hat I_0^{\prime\prime},
\end{eqnarray}
where $``\textrm{Int}\,A"$ and $``\textrm{co}\,A"$ are respectively the interior and the convex hull of a set $A$,  $B_{\textrm{aff}\,\mathcal K^{\hat\psi}}(\delta_0)$ is the closed ball in subspace $\textrm{aff}\,\mathcal K^{\hat\psi}$ with center $0\in\textrm{aff}\,\mathcal K^{\hat\psi}$ and radius $\delta_0$, and $D^2\hat\psi(\bar e)(y)=\big(D^2\hat\psi_1(\bar e)(y),\cdots,D^2\hat\psi_k(\bar e)(y)\big)^\top$.

\end{lem}

\textit{Proof.}\quad
First, since $E$ is convex, one can check by definition that  $T_E^{\flat(2)}(\bar e,y)$ is convex, and consequently $\mathcal K$ and $\mathcal K^{\hat\psi}$ are convex.

Second, we claim that
\begin{align}\label{cs1}
0\in \textrm{ri} \,\mathcal K^{\hat\psi},
\end{align}
where $ \textrm{ri} \,\mathcal K^{\hat\psi}$ is the interior of set $\mathcal K^{\hat\psi}$ relative to its affine hull (see \cite[p.44]{rock}). Consequently, $\textrm{aff}\,\mathcal K^{\hat\psi}$ is a subspace of $\mathbb R^k$.
By contradiction, we assume (\ref{cs1}) were not true.  Since the affine hull of $\mathcal K^{\hat\psi}$ is closed (see \cite[p. 44]{rock}), and ri $\mathcal K^{\hat\psi}$  is not empty (by \cite[Theorem 6.2, p. 45]{rock}),  it follows from \cite[Lemma 3.1]{bern} that, there exists $\xi\in\mathbb R^k\setminus\{0\}$ such that
\begin{align}\label{ccs10}
\xi^\top( D\hat\psi(\bar e)( x)+\frac{1}{2}D^2\hat\psi(\bar e)(y))\leq 0,\quad\forall x\in T_E^{\flat(2)}(\bar e,y).
\end{align}
Consequently we have
\begin{align}\label{ccs11}
(0,\xi^\top)(D\hat\Phi_{\hat I_0^{\prime\prime}}(\bar e)(x)+\frac{1}{2}D^2\hat\Phi_{\hat I_0^{\prime\prime}}(\bar e)(y))\leq 0,\quad\forall x\in T_E^{\flat(2)}(\bar e, y),
\end{align}
which contradicts the condition that $\mathcal K$ and $Z\times \{0\}$ are not separated by any linear functional. 

Third, if $D(\mathcal K^{\hat\psi})\leq 0$, then $D(\mathcal K^{\hat\psi})=0$, due to $\mathcal K^{\hat\psi}\neq\emptyset$. Consequently, we have $\mathcal K^{\hat\psi}=\{0\}$. Then, for any $\beta\in\mathbb R^k\setminus\{0\}$, (\ref{ccs10}) holds with $\xi$ replaced by $\beta$, and consequently (\ref{ccs11}) holds with $(0,\xi^\top)$ replaced by $(0,\beta^\top)$. A contradiction follows.

Finally, there exist $x_1,\cdots, x_{D(\mathcal K^{\hat\psi})+1}\in T_E^{\flat(2)}(\bar e,y)$ such that
\begin{equation}\label{c122}
0\in Int \,co\{D\hat\psi(\bar e)(x_l)+\frac{1}{2}D^2\hat\psi(\bar e)(y)\}_{l=1}^{D(\mathcal K^{\hat\psi})+1}.
\end{equation}

According to \cite[Lemma 3.1]{bern}, $(Z\times \{0\})\cap \mathcal K\neq \emptyset$. Then, there exist $z_0,z_1,\cdots,z_j\in(-\infty,0)$, $\lambda>0$ and $\tilde x\in T_E^{\flat(2)}(\bar e,y)$ such that
\begin{equation}\label{c12}\begin{array}{ll}
z_i=D\hat\phi_i(\bar e)(\tilde x)+\frac{1}{2}D^2\hat\phi_i(\bar e)(y),\quad&\textrm{if}\; \hat\phi_i(\bar e)=0,\,D\hat\phi_i(\bar e)(y)= 0;
\\ z_i-\lambda \hat\phi_i(\bar e)=0,&\textrm{if}\;\hat\phi_i(\bar e)<0;
\\ z_i-\lambda D\hat\phi_i(\bar e)(y)=0,&\textrm{if}\;\hat\phi_i(\bar e)=0,\,D\hat\phi_i(\bar e)(y)<0;
\end{array}\end{equation}
and
$
D\hat\psi(\bar e)(\tilde x)+\frac{1}{2}D^2\hat\psi(\bar e)(y)=0.
$
Let $\eta\in(0,1)$ be such that
$$\begin{array}{lll}
0&>&(1-\eta)\Big(D\hat\phi_i(\bar e)(x_l)+\frac{1}{2}D^2\hat\phi_i(\bar e)(y)\Big)+\eta\Big(D\hat\phi_i(\bar e)(\tilde x)+\frac{1}{2}D^2\hat\phi_i(\bar e)(y)\Big)
\\&=&D\hat\phi_i(\bar e)((1-\eta)x_l+\eta\tilde x)+\frac{1}{2}D^2\hat\phi_i(\bar e)(y),
\end{array}$$
for all $l=1,\cdots,D(\mathcal K^{\hat\psi})+1$ and  $i\in\{0,1,\cdots,j\}$ satisfying $\hat\phi_i(\bar e)=0$ and $D\hat\phi_i(\bar e)(y)=0$.
Set
$
h_l=\eta\tilde x+(1-\eta)x_l,
$ for $l=1,\cdots, k+1$.
Then (\ref{c200}) follows,
and (\ref{c119}) follows from (\ref{c122}).
 $\Box$
 
 The following results are respectively the first and second order necessary conditions of an minimizer of $(OP)$, and the idea of proving them partially comes from \cite[Theorem 4.1]{gb}.
 \begin{thm}\label{c223}
 Assume that $(C5)$  holds, and that $\bar e\in E$ is a solution to problem $(OP)$. Then, there exists $(\ell_{\hat\phi_0},\cdots,\ell_{\hat\phi_j},\ell_{\hat\psi}^\top)^\top\in\mathbb R^{1+j+k}\setminus\{0\}$ such that
 \begin{align}\label{lm}
\ell_{\hat\phi_i}\in(-\infty,0],\;i=0,\cdots,j;
\quad \ell_{\hat\phi_i}=0,\;\textrm{if}\;i\in \hat I_N;\quad \sum_{i=0}^j\ell_{\hat\phi_i}\hat\phi_i(\bar e)=0;
\\[2mm] \label{c400}
\sum_{i\in \hat I_
{A}}\ell_{\hat\phi_i}D\hat\phi_i(\bar e)(x)+\ell_{\hat\psi}^\top D\hat\psi(\bar e)(x)\leq 0,\;\forall x\in T_E^\flat(\bar e).
\end{align}
 \end{thm}
 
\begin{thm}\label{c126} Assume that $(C5)$ holds,  that $\bar e\in E$ is a solution to problem $(OP)$ with $\hat\phi_0(\bar e)=0$, and that  $y\in T_E^\flat(\bar e)$  satisfies  (\ref{c117}). Then, there exists $(\ell_{\hat\phi_0},\ell_{\hat\phi_1},\cdots,\ell_{\hat\phi_j},\ell_{\hat\psi}^\top)^\top\in \mathbb R^{1+j+k}\setminus\{0\}$ satisfying (\ref{lm}), (\ref{c400}),
\begin{eqnarray} \label{c2666}
\ell_{\hat\phi_i}=0,\;\textrm{if}\;i\notin \hat I_0^{\prime\prime},
\end{eqnarray}
and 
\begin{equation}\label{c141}\begin{array}{l}
\sum_{i\in \hat I_0^{\prime\prime}}\ell_{\hat\phi_i}D\hat\phi_i(\bar e)(x)+\ell_{\hat\psi}^\top D\hat\psi(\bar e)(x)
 +\sum_{i\in \hat I_0^{\prime\prime}}\frac{1}{2}\ell_{\hat\phi_i}D^2\hat\phi_i(\bar e)(y)
 \\[2mm]+\frac{1}{2}\ell_{\hat\psi}^\top D^2\hat\psi(\bar e)(y)\leq 0,\quad\forall\, x\in T_E^{\flat(2)}(\bar e, y).
\end{array}\end{equation}
\end{thm}

\begin{rem}\label{c225}
If $y\in T_E^\flat(\bar e)$ satisfies  (\ref{c117}), it is easy to see that the first order necessary condition  becomes trivial along the direction $y$: for any $\ell=(\ell_{\hat\phi_0},\cdots,\ell_{\hat\phi_j},\ell_{\hat\psi}^\top)^\top\in\mathbb R^{1+j+k}\setminus\{0\}$ satisfying  (\ref{lm}) and (\ref{c400}), it holds that $\sum_{i\in\hat I_A}\ell_{\hat\phi_i}D\hat\phi_i(\bar e)(y)+\ell_{\hat\psi}^\top D\hat\psi(\bar e)(y)=0$.
Thus Theorem \ref{c126} gives  further information of $\bar e$ along direction $y$. When $\bar e\in\textrm{Int}\,E$,  $0\in T_E^{\flat(2)}(\bar e,y)$, and consequently Theorem \ref{c126} is consistent with
 \cite[Theorem 4.1]{gb}. When $\bar e$ is on the boundary of $E$, $0\in T_E^{\flat(2)}(\bar e,y)$ is not always true, thus,  compared to  \cite[Theorem 4.1]{gb},  the first two terms of the left hand side of (\ref{c141}) are extra terms. 
\end{rem}

Since the proof of Theorem \ref{c223}
is analogous to that of Theorem \ref{c126}, we only prove Theorem \ref{c126} and give the key point of proving Theorem \ref{c223}:  The set
$\{D\hat\Phi_{\hat I_{A}}(\bar e)(x)|x\in T_E^\flat(\bar e)\}$ is separated from $\left((-\infty,0)^{j+1}-\{\lambda(\hat\phi(\bar e);\lambda>0\}\right)\times\{0\}$.

\medskip
\textbf{Proof of Theorem \ref{c126}}

{\it Step 1.}  We shall prove the case that $k>0$.

First, we claim that there exists $\ell\equiv(\ell_{\hat\phi_0},\cdots,\ell_{\hat\phi_j},\ell_{\hat\psi}^\top)^\top\in\mathbb R^{1+j+k}\setminus \{0\}$ such that (\ref{c129}) holds.

By contradiction, it follows from Lemma \ref{c115}  that (\ref{c119}) and (\ref{c200})  hold. 
Fix  $h_l$ ($l=1,\cdots,D(\mathcal K^{\hat\psi})+1$). Recall the definition of the second-order adjacent set (see Section \ref{na}). For any $\epsilon\to 0^+$, there exists $h^\epsilon_l\to h_l$ as $\epsilon\to 0^+$ such that $\bar e+\epsilon y+\epsilon^2 h_l^\epsilon\in E$. Then, there exists $\epsilon_0>0$ such that
$\bar e+\epsilon y+\epsilon^2h_l^\epsilon\in E$ for all $\epsilon\in[0,\epsilon_0]$ and $l=1,\cdots,D(\mathcal K^{\hat\psi})+1$. By the convexity of $E$, for any $x=\sum_{l=1}^{D(\mathcal K^{\hat\psi})+1}\nu_l h_l\in co\{h_1,\cdots,h_{D(\mathcal K^{\hat\psi})+1}\}$ with $(\nu_1,\cdots,\nu_{D(\mathcal K^{\hat\psi})+1})$ satisfying
\begin{equation}\label{c124}\sum_{l=1}^{D(\mathcal K^{\hat\psi})+1}\nu_l=1; \quad \nu_l\geq 0\quad\textrm{for}\; l=1,\cdots,D(\mathcal K^{\hat\psi})+1,\end{equation} 
it holds that
$\bar e+\epsilon y+\epsilon^2\sum_{l=1}^{D(\mathcal K^{\hat\psi})+1}\nu_lh_l^\epsilon \in E,$ for all $\epsilon\in[0,\epsilon_0]$.

Relation (\ref{c117}) and $(C5)$ imply that  exists $\epsilon_1\in(0,\epsilon_0]$ such that 
$$\begin{array}{ll}
&\Big|\epsilon^{-2}\hat\psi(\bar e+\epsilon y+\epsilon^2 \sum_{l=1}^{D(\mathcal K^{\hat\psi})+1}\nu_l h_l^\epsilon)-D\hat\psi(\bar e)(\sum_{l=1}^{D(\mathcal K^{\hat\psi})+1}\nu_lh_l)-\frac{1}{2}D^2\hat\psi(\bar e)(y)\Big|
\\\leq &\Big|\sum_{l=1}^{D(\mathcal K^{\hat\psi})+1}\nu_lD\hat\psi(\bar e)(h_l^\epsilon-h_l)\Big|
\\&+\Big|\epsilon^{-2}\hat\psi(\bar e+\epsilon y+\epsilon^2 \sum_{l=1}^{D(\mathcal K^{\hat\psi})+1}\nu_l h_l^\epsilon)-D\hat\psi(\bar e)(\sum_{l=1}^{D(\mathcal K^{\hat\psi})+1}\nu_lh_l^\epsilon)-\frac{1}{2}D^2\hat\psi(\bar e)(y)\Big|
\\<&\delta_0  ,\;\forall\epsilon\in[0,\epsilon_1],
\end{array}$$
where $\nu_1,\cdots,\nu_{D(\mathcal K^{\hat\psi})+1}$ satisfy 
(\ref{c124}).

According to (\ref{c117}), (\ref{c200}) and $(C5)$, one can find $\epsilon_2\in(0,\epsilon_1]$ such that, for any $\epsilon\in[0,\epsilon_2]$, the following relations hold: for $i\in\hat I_0^{\prime\prime}$,
\begin{equation}\label{c220}
\begin{array}{ll}
&\hat\phi_i(\bar e+\epsilon y+\epsilon^2\sum_{l=1}^{D(\mathcal K^{\hat\psi})+1}\nu_l h_l^\epsilon)
\\=&\epsilon^2\sum_{l=1}^{D(\mathcal K^{\hat\psi})+1}\nu_l\left(D\hat\phi_i(\bar e)(h_l)+\frac{1}{2}D^2\hat\phi_i(\bar e)(y)
\right)
\\&+\epsilon^2\sum_{l=1}^{D(\mathcal K^{\hat\psi})+1}\nu_lD\hat\phi_i(\bar e)(h_l^\epsilon-h_l)+\Big[\hat\phi_i(\bar e+\epsilon y+\epsilon^2\sum_{l=1}^{D(\mathcal K^{\hat\psi})+1}\nu_l h_l^\epsilon)
\\&-\epsilon^2\sum_{l=1}^{D(\mathcal K^{\hat\psi})+1}\nu_l\left(D\hat\phi_i(\bar e)(h_l^\epsilon)+\frac{1}{2}D^2\hat\phi_i(\bar e)(y)
\right)\Big]
\\<&0;
\end{array}\end{equation}
 for $i\notin \hat I_0^{\prime\prime}$,
\begin{align}\label{c221}\begin{array}{ll}
 &\hat\phi_i(\bar e+\epsilon y+\epsilon^2\sum_{l=1}^{D(\mathcal K^{\hat\psi})+1}\nu_l h_l^\epsilon)
\\=&\hat\phi_i(\bar e)+\epsilon D\hat\phi_i(\bar e)(y)+\epsilon^2\sum_{l=1}^{D(K^{\hat\psi})+1}\nu_l\Big(D\hat\phi_i(\bar e)h_l^\epsilon+\frac{1}{2}D^2\hat\phi_i(\bar e)(y)\Big)
\\&+\Big[\hat\phi_i(\bar e+\epsilon y+\epsilon^2\sum_{l=1}^{D(\mathcal K^{\hat\psi})+1}\nu_l h_l^\epsilon)
-\hat\phi_i(\bar e)-\epsilon D\hat\phi_i(\bar e)(y)
\\&-\epsilon^2\sum_{l=1}^{D(K^{\hat\psi})+1}\nu_l\Big(D\hat\phi_i(\bar e)h_l^\epsilon+\frac{1}{2}D^2\hat\phi_i(\bar e)(y)\Big)\Big]
\\<&0.
\end{array}\end{align}

Then, from the above relations and (\ref{c119}), we can define a map 
\begin{align*}
G: co\{D\hat\psi(\bar e)(h_l)+\frac{1}{2}D^2\hat\psi(\bar e)(y)\}_{l=1}^{D(\mathcal K^{\hat\psi})+1}\to 
\\  co\{D\hat\psi(\bar e)(h_l)+\frac{1}{2}D^2\hat\psi(\bar e)(y)\}_{l=1}^{D(\mathcal K^{\hat\psi})+1}
\end{align*}
 by
$$\begin{array}{ll}
&G\Big(D\hat\psi(\bar e)(\sum_{l=1}^{D(\mathcal K^{\hat\psi})+1}\nu_lh_l)+\frac{1}{2}D^2\hat\psi(\bar e)(y)\Big)
\\=&-\epsilon_2^{-2}\hat\psi(\bar e+\epsilon_2 y+\epsilon_2^2 \sum_{l=1}^{D(\mathcal K^{\hat\psi})+1}\nu_l h_l^{\epsilon_2})+D\hat\psi(\bar e)(\sum_{l=1}^{D(\mathcal K^{\hat\psi})+1}\nu_lh_l)+\frac{1}{2}D^2\hat\psi(\bar e)(y),
\end{array}$$
for all $\nu_1,\cdots,\nu_{D(\mathcal K^{\hat\psi})+1}$ satisfying (\ref{c124}). Obviously $G$ is continuous and 
$ co\{D\hat\psi(\bar e)(h_l)+\frac{1}{2}D^2\hat\psi(\bar e)(y)\}_{l=1}^{D(\mathcal K^{\hat\psi})+1}$ is convex and compact. By Brouwer fixed point theorem, there exists $\nu_1^\star,\cdots, \nu_{D(\mathcal K^{\hat\psi})+1}^\star$ satisfying (\ref{c124}) such that
$
G\Big(D\hat\psi(\bar e)\Big(\sum_{l=1}^{D(\mathcal K^{\hat\psi})+1}\nu_l^\star h_l\Big)+\frac{1}{2}D^2\hat\psi(\bar e)(y)\Big)
=D\hat\psi(\bar e)(\sum_{l=1}^{D(\mathcal K^{\hat\psi})+1}\nu_l^\star h_l)+\frac{1}{2}D^2\hat\psi(\bar e)(y),
$
which implies
$\hat\psi(\bar e+\epsilon_2 y+\epsilon_2^2 \sum_{l=1}^{D(\mathcal K^{\hat\psi})+1}\nu_l^\star h_l^{\epsilon_2})=0$.
Recalling (\ref{c220}) and  (\ref{c221}), we obtain that $\bar e+\epsilon_2 y+\epsilon_2^2 \sum_{l=1}^{D(\mathcal K^{\hat\psi})+1}\nu_l^\star h_l^{\epsilon_2}$ satisfies (\ref{c5}) and $\hat\phi_0(\bar e+\epsilon_2 y+\epsilon_2^2 \sum_{l=1}^{D(\mathcal K^{\hat\psi})+1}\nu_l^\star h_l^{\epsilon_2})<0$, which contradicts the optimality of $\bar e$.

Second, from (\ref{c129}) and the special structure of (\ref{c201}), we obtain (\ref{lm}), (\ref{c2666}), and 
$$
\ell^\top (D\hat\Phi_{\hat I_0^{\prime\prime}}(\bar e)(x)+\frac{1}{2}D^2\hat\Phi_{\hat I_0^{\prime\prime}}(\bar e)(y))\leq\inf_{z\in Z} \ell^\top(z^\top,0)^\top=0,\;\forall x\in T_E^{\flat(2)}(\bar e,y),
$$
which implies (\ref{c141}).

Finally, from \cite[Proposition 4.2.1]{AF} and \cite[Lemma 2.4]{fdt} we have
\begin{align}\label{cs4}
T^{\flat(2)}_E(\bar e,y)=T^{\flat(2)}_E(\bar e,y)+T^\flat_E(\bar e).
\end{align} If (\ref{c400}) were not true, there would exist $x_0\in T_E^\flat(\bar e)$ such that
\\$\sum_{i\in \hat I_{A}}\ell_{\hat\phi_i}D\hat\phi_i(\bar e)(x_0)+\ell_{\hat\psi}^\top D\hat\psi(\bar e)(x_0)>0$. Fix any $\sigma \in T_E^2(\bar e,y)$. Choosing $\lambda>0$ big enough such that 
$$\begin{array}{l}
\sum_{i\in \hat I_0^{\prime\prime}}\ell_{\hat \phi_i}D\hat\phi_i(\bar e)(\sigma+\lambda x_0)+\ell_{\hat\psi}^\top D\hat\psi(\bar e)(\sigma+\lambda x_0)
+\sum_{i\in\hat I_0^{\prime\prime}}\frac{1}{2}\ell_{\hat\phi_i}D^2\hat\phi_i(\bar e)(y)
 \\+\frac{1}{2}\ell_{\hat\psi}^\top D^2\hat\psi(\bar e)(y)>0,
\end{array}$$
which contradicts (\ref{c141}),   
and the proof is concluded.

{\it Step 2.}  For the case $k=0$, there exists $(\ell_{\hat\phi_0},\cdots,\ell_{\hat\phi_j})^\top\in\mathbb R^{1+j}\setminus\{0\}$ such that
\begin{align*}
\sum_{l=0}^{j}\ell_{\hat\phi_l}\beta_l\leq \sum_{l=0}^j\ell_{\hat\phi_l}z_l,\quad\forall\,(\beta_0,\cdots,\beta_j)^\top\in\mathcal K,\;(z_0,\cdots,z_j)^\top\in Z.
\end{align*}
If it were not true, by \cite[Lemma 3.1]{bern} we would have $\mathcal K\cap Z\neq\emptyset$. Then, there exists $\tilde x\in T_E^{\flat(2)}(\bar e,y)$, $\lambda>0$ and $(z_0,\cdots,z_j)^\top\in(-\infty,0)^{1+j}$ such that (\ref{c12}) holds. Thus, for any $i=0,\cdots,j$, we obtain from $(C5)$ and (\ref{c12}) that, there exists $\tilde\epsilon>0$ such that, for any $\epsilon\in[0,\tilde\epsilon]$ the following relation 
\begin{align*}
&\hat\phi_i(\bar e+\epsilon y+\epsilon^2\tilde x)
\\=&\hat \phi_i(\bar e)+\epsilon D\hat\phi_i(\bar e)(y)+\epsilon^2 D\hat \phi_i(\bar e)(\tilde x)+\frac{\epsilon^2}{2}D^2\hat\phi_i(\bar e)(y)+\Big[\hat\phi_i(\bar e+\epsilon y+\epsilon^2\tilde x)
\\&-\hat \phi_i(\bar e)-\epsilon D\hat\phi_i(\bar e)(y)-\epsilon^2 D\hat \phi_i(\bar e)(\tilde x)-\frac{\epsilon^2}{2}D^2\hat\phi_i(\bar e)(y)\Big]
\\<&0,
\end{align*}
holds for $i=0,1,\cdots,j$, which contradicts the optimality of $\bar e$.
$\Box$

\subsection{Proof of Theorem \ref{c230} }\label{psoc}

We need the following  lemmas.

\begin{lem}\label{l10}
Assume $U\subset \mathbb R^m$ is closed. Fix $\bar u(\cdot)\in \mathcal U$. Let $v(\cdot)\in L^1(0,T;\mathbb R^m)$ be such that $v(t)\in T_U^\flat(\bar u(t))$ a.e. $t\in[0,T]$. Assume 
 there exist a positive constant  $\epsilon_0$ and $\ell(\cdot)\in L^i(0,T;\mathbb R^m)$ ($i=1$ or $2$) such that (\ref{l4}) holds.
Fix any $\sigma(\cdot)\in L^1(0,T;\mathbb R^m)$ such that $\sigma(t)\in T_U^{\flat(2)}(\bar u(t), v(t))$ a.e. $t\in[0,T]$. Then, for any $\epsilon\in(0,\epsilon_0]$, there exists $\sigma_\epsilon(\cdot)\in L^i(0,T;\mathbb R^m)$ such that
$u_\epsilon(t):=\bar u(t)+\epsilon v(t)+\epsilon^2\sigma_\epsilon(t)\in U$ and  $\lim_{\epsilon\to 0^+}\sigma_\epsilon(t)=\sigma(t)$ a.e. $t\in[0,T]$, and 
\begin{equation}\label{l7}
\|\sigma_\epsilon(\cdot)\|_{L^i(0,T;\mathbb R^m)}\leq \|\ell\|_{L^i(0,T;\mathbb R^m)}+2\|\sigma\|_{L^i(0,T;\mathbb R^m)}.
\end{equation}
\end{lem}

\textbf{Proof}\quad
Since $U$ is closed, by \cite[Corollary 8.2.13, p. 317]{AF}, for every $\epsilon>0$, there exist  measurable  functions $\omega_\epsilon, z_\epsilon: [0,T]\to U$ such that 
\begin{align}
\label{l2}dist_U(\bar u(t)+\epsilon v(t)+\epsilon^2 \sigma(t))=|\bar u(t)+\epsilon v(t)+\epsilon^2 \sigma(t)-\omega_\epsilon(t)|\quad a.e.\,t\in[0,T],
\\
\label{l3}
a_\epsilon(t):=dist_U(\bar u(t)+\epsilon v(t))=|\bar u(t)+\epsilon v(t)-z_\epsilon(t)|\quad a.e.\,t\in [0,T].
\end{align}
Set
$
\sigma_\epsilon(t)=\frac{1}{\epsilon^2}\big(\omega_\epsilon(t)-\bar u(t)-\epsilon v(t)\big)$ for $ t\in[0,T]$.
Then, we have $\bar u(t)+\epsilon v(t)+\epsilon^2\sigma_\epsilon(t)\in U$ a.e   $t\in [0,T]$.
Since $\sigma(t)\in T_U^{\flat(2)}(\bar u(t),v(t))$ a.e. $t\in [0,T]$, recalling (\ref{l2}) and (\ref{ccs90}),  we have
\begin{equation}\label{l6}
0=\lim_{\epsilon\to 0^+}\frac{1}{\epsilon^2}|\bar u(t)+\epsilon v(t)+\epsilon^2 \sigma(t)-\omega_\epsilon(t)|=\lim_{\epsilon\to 0^+}|\sigma_\epsilon(t)-\sigma(t)|,\quad a.e.\,t\in[0,T].
\end{equation}
Applying (\ref{l2}), (\ref{l3}) and (\ref{l4}), we have, for $\epsilon\in(0,\epsilon_0]$,
$$\begin{array}{l}
\epsilon^2(|\sigma_\epsilon(t)|-|\sigma(t)|)\leq |\epsilon^2\sigma(t)-\epsilon^2\sigma_\epsilon(t)|= |\bar u(t)+\epsilon v(t)+\epsilon^2\sigma(t)-\omega_\epsilon(t)|
\\ \leq |\bar u(t)+\epsilon v(t)+\epsilon^2\sigma(t)-z_\epsilon(t)|
\leq a_\epsilon(t)+\epsilon^2|\sigma(t)|\leq \epsilon^2(\ell(t)+|\sigma(t)|),\quad a.e.\,t\in[0,T],
\end{array}$$
which implies that, $\sigma_\epsilon(\cdot)\in L^i(0,T;\mathbb R^m)$ and (\ref{l7}), if $\ell(\cdot)\in L^i(0,T;\mathbb R^m)$ ($i=1,2$).
$\Box$

\begin{lem}\label{ccs91}
Assume $U\subset\mathbb R^m$ is closed. Fix $\bar u(\cdot)\in L^2(0,T;\mathbb R^m)\cap \mathcal U$. Let $v(\cdot)\in L^2(0,T;\mathbb R^m) $ be such that $v(t)\in T^\flat_U(\bar u(t))$ a.e. in $[0,T]$. Then, $v(\cdot)\in T^\flat_{L^2(0,T;\mathbb R^m)\cap\mathcal U}(\bar u(\cdot))$. Moreover, for any $\sigma(\cdot)\in L^2(0,T;\mathbb R^m)$ such that $\sigma(t)\in T^{\flat(2)}_U(\bar u(t),v(t))$ a.e. in $[0,T]$ and (\ref{l4}) holds for some $\ell(\cdot)\in L^2(0,T;\mathbb R^m)$ and $\epsilon_0>0$, it holds that $\sigma(\cdot)\in T^{\flat(2)}_{\mathcal U\cap L^2(0,T;\mathbb R^m)}(\bar u(\cdot),\\ v(\cdot))$.
\end{lem}

{\it Proof.}\; First, by \cite[Corollary 8.2.13, p. 317]{AF}, for each $\epsilon>0$, there exists a measurable map $v_\epsilon:[0,T]\to \mathbb R^m$ such that $v_\epsilon(t)\in \frac{U-\{\bar u(t)\}}{\epsilon}$ and $dist_{\frac{U-\{\bar u(t)\}}{\epsilon}}v(t)=|v(t)-v_\epsilon(t)|$ for almost every $t\in[0,T]$. It follows from (\ref{ccs94}) that $\lim_{\epsilon\to 0^+}|v(t)-v_\epsilon(t)|=0$ a.e.  $t\in[0,T]$.  Then, we have
\begin{align*}
|v_\epsilon(t)|-|v(t)|\leq |v_\epsilon(t)-v(t)|\leq |v(t)-\frac{1}{\epsilon}(\bar u(t)-\bar u(t))|=|v(t)|.
\end{align*}
We obtain from Lebesgue’s
dominated convergence theorem that $\lim_{\epsilon\to 0^+} v_\epsilon(\cdot)=v(\cdot)$ in $L^2(0,T;\mathbb R^m)$, consequently we have $v(\cdot)\in T^\flat_{\mathcal U\cap L^2(0,T;\mathbb R^m)}(\bar u(\cdot))$.

Then, it follows from Lemma \ref{l10} that, for any $\epsilon>0$, there exists $\sigma_\epsilon(\cdot)\in L^2(0,T;\mathbb R^m)$ such that $\lim_{\epsilon\to 0^+}\sigma_\epsilon(\cdot)=\sigma(\cdot)$ in $L^2(0,T;\mathbb R^m)$ and $\bar u(\cdot)+\epsilon v(\cdot)+\epsilon^2 \sigma_\epsilon(\cdot)\in \mathcal U\cap L^2(0,T;\mathbb R^m)$, and consequently $\sigma(\cdot)\in  T^{\flat(2)}_{\mathcal U\cap L^2(0,T;\mathbb R^m)}(\bar u(\cdot), v(\cdot))$. The proof is concluded. 
$\Box$

\medskip

 Then, we are going to prove Theorem  \ref{c230}.

\textbf{Proof of Theorem \ref{c230}}\quad
First, we shall transform problem $(OCP)$ to an optimization problem. Assume $(\bar u(\cdot), \bar y(\cdot))$ is an optimal pair for problem $(OCP)$ with $\bar u(\cdot)\in L^2(0,T;\mathbb R^m)$.  For any $(Y,u(\cdot))\in T_{\bar y(0)}M\times L^2(0,T;\mathbb R^m)$, set
\begin{align}\label{ccs110}\begin{array}{l}
\hat\phi_i(Y,u(\cdot))\equiv \phi_i\left(y_u(0;exp_{\bar y(0)}Y),y_u(T;exp_{\bar y(0)}Y)\right),\;i=1,\cdots,j,
\\ \hat\phi_0(Y,u(\cdot))\equiv \phi_0\left(y_u(0;exp_{\bar y(0)}Y),y_u(T;exp_{\bar y(0)}Y)\right)-\phi_0(\bar y(0),\bar y(T)),
\\ \hat\psi(Y,u(\cdot))\equiv\psi\left(y_u(0;exp_{\bar y(0)}Y),y_u(T;exp_{\bar y(0)}Y)\right),
\end{array}
\end{align}
where   $y_u(\cdot;x)$ is the solution to (\ref{s2}) with initial state $x\in M$ and control $u(\cdot)$, and $exp_x\cdot$ is the exponential map at $x$ (see Section \cite[Section 2.1]{cdz}).

We obtain  from the optimality of $(\bar u(\cdot),\bar y(\cdot))$ for problem $(OCP)$ that, $(0,\bar u(\cdot))\in T_{\bar y(0)}M\times \mathcal U$ is the solution to the following optimization problem
\begin{description}
\item[$(\widetilde{OCP})$] Find $(Y,u(\cdot))\in T_{\bar y(0)}M\times\big( L^2(0,T;\mathbb R^m)\cap\mathcal U\big)$ minimizes $\hat\phi_0(Y,u(\cdot))$
subject to  $\hat\phi_i(Y,u(\cdot))\leq 0$ for $i=1,\cdots,j$,  $ \hat\psi(Y,u(\cdot))=0$ and  $ (Y,u(\cdot))\in  T_{\bar y(0)}M\times( \mathcal U\cap L^2(0,T;\mathbb R^m))$.

\end{description}

Second, we shall check that condition $(C5)$ holds. Fix  $(V,v(\cdot))\in T_{\bar y(0)}M\times L^2(0,T;\mathbb R^m)$. For $\epsilon>0$,  we denote by $y(\cdot;\exp_{\bar y(0)}\epsilon V, \bar u(\cdot)+\epsilon v(\cdot))$ the solution to    (\ref{s2}) corresponding to the initial state $\exp_{\bar y(0)}\epsilon V$ and the control $\bar u(\cdot)+\epsilon v(\cdot)$. For $i=0,1,\cdots,j$, 
we obtain from Proposition \ref{439} that
\begin{align*}
&\hat\phi_i(\epsilon V,\bar u(\cdot)+\epsilon v(\cdot))-\hat\phi_i(0,\bar u(\cdot))
\\=&\epsilon\Big[\nabla_1\phi_i(\bar y(0),\bar y(T))(V)+\nabla_2\phi_i(\bar y(0),\bar y(T))(X_{v,V}(T))\Big]+o(\epsilon),
\end{align*}
where $X_{v,V}(\cdot)$ is the solution to (\ref{c227}) with $X_{v,V}(0)=V$.
This implies that $\hat\phi_i$ is Fr$\acute {\textrm e}$chet differentiable at $(0,\bar u(\cdot))$, and its Fr$\acute {\textrm e}$chet derivative is as follows
\begin{align}\label{ccs111}
D\hat\phi_i(0,\bar u(\cdot))(V,v(\cdot))=\nabla_1\phi_i(\bar y(0),\bar y(T))(V)+\nabla_2\phi_i(\bar y(0),\bar y(T))(X_{v,V}(T)).
\end{align}
Similarly we can show that $\hat\psi$ is  Fr$\acute {\textrm e}$chet differentiable at $(0,\bar u(\cdot))$, and its Fr$\acute {\textrm e}$chet derivative is given by
\begin{align}\label{ccs112}
D\hat\psi(0,\bar u(\cdot))(V,v(\cdot))=\nabla_1\psi(\bar y(0),\bar y(T))(V)+\nabla_2\psi(\bar y(0),\bar y(T))(X_{v,V}(T)),
\end{align}
where $\nabla_i\psi$ ($i=1,2$) is defined in (\ref{pcd}). 

Fix any $(W,\sigma(\cdot))\in T_{\bar y(0)}M\times L^2(0,T;\mathbb R^m)$. Denote by $X_{\sigma, W}(\cdot)$ the solution to (\ref{c227}) with $v(\cdot)$ replaced by $\sigma(\cdot)$ and $X_{\sigma,W}(0)=W$. For any $\epsilon>0$,  we denote by $y_\epsilon(\cdot)$ the solution to (\ref{s2}) with initial state $\exp_{\bar y(0)}(\epsilon V+\epsilon^2 W)$ and control $\bar u(\cdot)+\epsilon v(\cdot)+\epsilon^2 \sigma(\cdot)$. Denote by $Y_{\sigma W}^{X_{v,V}}(\cdot)$ the solution to (\ref{o2}) with $(\sigma_\epsilon(\cdot), X_v(\cdot))$ replaced by $(\sigma(\cdot), X_{v,V}(\cdot))$. We employ the notations $V_\epsilon(\cdot)$ and $\beta(\cdot; t)$ ($t\in[0,T]$) given respectively by (\ref{438}) and (\ref{ccs100}).   Note that (\ref{ccs34}) still holds. Fix $\alpha>0$. It follows from  Proposition \ref{439} that, there exists $\epsilon_1>0$ such that
\begin{align}\label{ccs102}
V_\epsilon(t)=\epsilon X_{v,V}(t)+\epsilon^2 Y_{\sigma W}^{X_{v,V}}(t)+\gamma_\epsilon(t),\;\forall\,t\in[0,T],\;\epsilon\in[0,\epsilon_1],
\end{align}
with
\begin{align}\label{ccs1111}
|\gamma_\epsilon(t)|\leq \frac{\alpha}{2K}\epsilon^2,\quad\forall\,\epsilon\in[0,\epsilon_1],
\end{align}
where  constant $K$ is given in condition $(C2)$.

Set by
\begin{align}\label{ccs114}\begin{array}{ll}
&D^2\hat\phi_i(0,\bar u(\cdot))(V,v(\cdot))
\\=&\nabla_1^2\phi_i(\bar y(0),\bar y(T))(V,V)+2\nabla_2\nabla_1\phi_i(\bar y(0),\bar y(T))(V,X_{v,V}(T))
\\&+\nabla_2^2\phi_i(\bar y(0),\bar y(T))(X_{v,V}(T),X_{v,V}(T))+2\nabla_2\phi_i(\bar y(0),\bar y(T))(Y_{00}^{X_{v,V}}(T)),
\end{array}\end{align}
where $Y_{00}^{X_{v,V}}(\cdot)$ is the solution to (\ref{o2}) with $\sigma(\cdot)=0$ and $W=0$, and with $X_v(\cdot)$ replaced by $X_{v,V}(\cdot)$.
It is easy to check that 
\begin{align}\label{ccs117}
Y_{\sigma W}^{X_{v,V}}(t)=Y_{00}^{X_{v,V}}(t)+X_{\sigma,W}(t),\quad\forall\,t\in[0,T],
\end{align}
where $X_{\sigma,W}(\cdot)$ is the solution to 
(\ref{c227}) with initial state $X_{\sigma,W}(0)=W$, and with $v(\cdot)$ replaced by $\sigma(\cdot)$.

 For $i=0,1, \cdots,j$, we obtain by Newton-Leibniz formula, exchange of integral variables, (\ref{ccs102}) and (\ref{ccs117}) that
\begin{align*}
&\hat\phi_i(\epsilon V+\epsilon^2W,\bar u(\cdot)+\epsilon v(\cdot)+\epsilon^2\sigma(\cdot))-\hat\phi_i(0,\bar u(\cdot))-\epsilon D\hat\phi_i(0,\bar u(\cdot))(V,v(\cdot))
\\&-\epsilon^2\Big[D\hat\phi_i(0,\bar u(\cdot))(W,\sigma(\cdot))+\frac{1}{2}D^2\hat\phi_i(0,\bar u(\cdot))(V,v(\cdot))\Big]=
L_i^\epsilon,
\end{align*}
where 
\begin{align*}
L_i^\epsilon=&\int_0^1\Big[\nabla_1^2\phi_i(
\beta(\tau;0),\beta(\tau;T))\Big(\frac{\partial}{\partial \tau}\beta(\tau;0),\frac{\partial}{\partial \tau}\beta(\tau;0)\Big)-\nabla_1^2\phi_i(\bar y(0),\bar y(T))(V,V)\epsilon^2
\\&+2\nabla_2\nabla_1\phi_i(
\beta(\tau;0),\beta(\tau;T))\Big(\frac{\partial}{\partial \tau}\beta(\tau;0),\frac{\partial}{\partial \tau}\beta(\tau;T)\Big)-2\epsilon^2\nabla_2
\nabla_1\phi_i(\bar y(0),\bar y(T))(V,
\\&X_{v,V}(T))+\nabla_2^2\phi_i
(
\beta(\tau;0),\beta(\tau;T))\Big(\frac{\partial}{\partial \tau}\beta(\tau;T),\frac{\partial}{\partial \tau}\beta(\tau;T)\Big)
\\&-\epsilon^2\nabla_2^2\phi_i(\bar y(0),\bar y(T))(X_{v,V}(T), X_{v,V}(T))\Big](1-\tau)d\tau+\nabla_2\phi_i(\bar y(0),\bar y(T))(\gamma_\epsilon(T)).
\end{align*}
By  (\ref{ccs34}) and (\ref{ccs102}), we  have
\begin{align*}
L_i^\epsilon=&\epsilon^2\int_0^1\Big[\nabla_1^2
\phi_i(\beta(\tau;0),\beta(\tau;T))(L_{\bar y(0)\beta(\tau;0)}V,
L_{\bar y(0)\beta(\tau;0)}V)-\nabla_1^2\phi_i(\bar y(0),\bar y(T))(V,V)
\\&+2\nabla_2\nabla_1\phi_i(\beta(\tau;0),\beta(\tau;T))(L_{\bar y(0)\beta(\tau;0)}V,L_{\bar y(T)\beta(\tau;T)}X_{v,V}(T))
\\&-2\nabla_2\nabla_1\phi_i(\bar y(0),\bar y(T))(V,X_{v,V}(T))+\nabla_2^2\phi_i(\beta(\tau;0),\beta(\tau;T))(L_{\bar y(T)\beta(\tau;T)}X_{v,V}(T),
\\&L_{\bar y(T)\beta(\tau;T)}X_{v,V}(T))-\nabla_2^2\phi_i(\bar y(0),\bar y(T))(X_{v,V}(T),X_{v,V}(T))\Big](1-\tau)d\tau+o(\epsilon^2)
\\&+\nabla_2\phi_i(\bar y(0),\bar y(T))(\gamma_\epsilon(T)).
\end{align*}
Applying Lebesgue’s
dominated convergence theorem, ($C2$), ($C3$),  (\ref{ccs111}) and (\ref{ccs117}) to the above identity, we obtain that, there exists $\tilde\epsilon_0\in(0,\epsilon_1]$ such that  $|L_i^\epsilon|\leq \alpha\epsilon^2$ for all $\epsilon\in[0,\tilde\epsilon_0]$. Similarly one can show that, there exists $\epsilon_0\in(0,\tilde \epsilon_0]$ such that
\begin{align*}
&\Big|\hat\psi(\epsilon V+\epsilon^2W,\bar u(\cdot)+\epsilon v(\cdot)+\epsilon^2\sigma(\cdot))-\hat\psi(0,\bar u(\cdot))-\epsilon D\hat\psi(0,\bar u(\cdot))(V,v(\cdot))
\\&-\epsilon^2[D\hat\psi(0,\bar u(\cdot))(W,\sigma(\cdot))+\frac{1}{2}D^2\hat\psi(0,\bar u(\cdot))(V,v(\cdot))\Big]\Big|\leq \alpha\epsilon^2,
\end{align*}
for all $\epsilon\in[0,\epsilon_0]$,
where
\begin{align}\label{ccs115}\begin{array}{ll}
&D^2\hat\psi(0,\bar u(\cdot))(V,v(\cdot))
\\=&\nabla_1^2\psi(\bar y(0),\bar y(T))(V,V)+2\nabla_2\nabla_1\psi(\bar y(0),\bar y(T))(V,X_{v,V}(T))
\\&+\nabla_2^2\psi(\bar y(0),\bar y(T))(X_{v,V}(T),X_{v,V}(T))+2\nabla_2\psi(\bar y(0),
\bar y(T))(Y_{00}^{X_{v,V}}(T)).
\end{array}\end{align}
Thus, for problem $(\widetilde{OCP})$, condition $(C5)$ holds.

Third, we shall use Theorem \ref{c115} to prove Theorem \ref{f}.
Recalling (\ref{c246}), (\ref{c226}), (\ref{c231}) and (\ref{ccs110}), we have $I_{A}= \{i\in\{1,\cdots,j\}|\,\hat\phi_i(0,\bar u(\cdot))=0\}\cup\{0\}$ and $I_N=\{i\in\{1,\cdots,j\}|\,\hat\phi_i(0,\bar u(\cdot))<0\}$.
Applying Theorem \ref{c223} to problem $(\widetilde{OCP})$, we  can find $\ell=(\ell_{\phi_0},\cdots,\ell_{\phi_j},\ell_\psi)\in\mathbb R^{1+j+k}\setminus\{0\}$ satisfying   (\ref{c260}) and
\begin{align}\label{c229}\begin{array}{l}
\sum_{i\in  I_{A}}\ell_{\phi_i}\left(\nabla_1\phi_i(\bar y(0),\bar y(T))(Y)+\nabla_2\phi_i(\bar y(0),\bar y(T))(X_{w,Y}(T))\right)
\\[2mm]+\ell_\psi^\top\left(\nabla_1\psi(\bar y(0),\bar y(T))(Y)+\nabla_2\psi(\bar y(0),\bar y(T))(X_{w,Y}(T))\right)\leq 0,
\end{array}
\end{align}
for all $(Y,w(\cdot))\in T_{\bar y(0)}M\times T_{\mathcal U\cap L^2(0,T;\mathbb R^m)}^\flat(\bar u(\cdot))$, where $X_{w,Y}(\cdot)$ is the solution to (\ref{c227}) with $X_{w,Y}(0)=Y$, and
we have used (\ref{ccs111}) and (\ref{ccs112}).  Let $p^\ell(\cdot)$ be the solution to (\ref{c243}).
Inserting (\ref{c243}) into (\ref{c229}) and integrating by parts, we can obtain from Lemma \ref{ccs91} that, $\int_0^T\nabla_u H[t,\ell](w(t))\leq 0$ holds for all $w(\cdot)\in L^2(0,T;\mathbb R^m)$ with $w(t)\in T_U^\flat(\bar u(t))$ a.e. $t\in[0,T]$, and that (\ref{icd1}) stands.  Applying  needle variation to this inequality, we obtain (\ref{c240}). Thus,  $\ell$ is a Lagrange multiplier in the sense of convex variation.

Finally, we shall employ Theorem \ref{c126} to prove Theorem \ref{c230}. Assume that $v(\cdot)\in L^2(0,T;\mathbb R^m)$ is a singular direction in the sense of convex variation, with $X_v(\cdot)$ satisfying (\ref{c266}) and (\ref{c227}), and that (\ref{l4}) holds for some  $\epsilon_0>0$ and $\ell(\cdot)\in L^2(0,T;\mathbb R^m)$.
Recall (\ref{ccs114}) and (\ref{ccs115}). It follows from Theorem \ref{c126} and Lemma \ref{ccs91} that,   there exist a Lagrange mulitplier in the sense of convex variation $\ell=(\ell_{\phi_0}, \ell_{\phi_1},\cdots, \ell_{\phi_j}, \ell_\psi^\top)^\top\in\mathbb R^{1+j+k}\setminus\{0\}$ satisfying 
(\ref{c232})  and 
\begin{equation}\label{c233}
\begin{array}{l}
\sum_{i\in  I_0^{\prime\prime}}\ell_{\phi_i}\Big(\nabla_1\phi_i(\bar y(0),\bar y(T))(W)+\nabla_2\phi_i(\bar y(0),\bar y(T))(X_{\sigma,W}(T)+Y_{00}^{X_{v}}(T))\Big)
\\+\ell_\psi^\top\Big(\nabla_1\psi(\bar y(0),\bar y(T))(W)+\nabla_2\psi(\bar y(0),\bar y(T))(X_{\sigma,W}(T)+Y_{00}^{X_{v}}(T))\Big)
\\+\frac{1}{2}\sum_{i\in  I_0^{\prime\prime}}\ell_{\phi_i}\Big( 
\nabla_1^2\phi_i(\bar y(0),\bar y(T))(X_v(0),X_v(0))+2\nabla_2\nabla_1\phi_i(\bar y(0),\bar y(T))(X_v(0),
\\X_{v}(T))+\nabla_2^2\phi_i(\bar y(0),\bar y(T))(X_{v}(T),X_{v}(T))\Big)+\frac{1}{2}\ell_\psi^\top\Big(\nabla_1^2\psi(\bar y(0),\bar y(T))(X_v(0),X_v(0))
\\ +2\nabla_2\nabla_1\psi(\bar y(0),\bar y(T))(X_v(0),X_{v}(T))+\nabla_2^2\psi(\bar y(0),\bar y(T))(X_{v}(T),X_{v}(T))\leq 0,
\end{array}
\end{equation}
for all $(W,\sigma)\in T_{\bar y(0)}M\times L^2(0,T;\mathbb R^m)$ with $\sigma(t)\in T_U^{\flat(2)}(\bar u(t),v(t))$ a.e. $t\in[0,T]$, where $Y_{00}^{X_v}(\cdot)$ is the solution to (\ref{o2}) with $\sigma(\cdot)=0$ and $W=0$.

Recall that $ p^\ell(\cdot)$ solves (\ref{c243}) with initial data (\ref{icd1}). We obtain from (\ref{c243}), (\ref{o2}), (\ref{ccs117}) and integration by parts over $[0,T]$ that
$$\begin{array}{lll}
0&\geq& - p^\ell(0)(W)+p^\ell(T)(Y_{\sigma W}^{X_{v,V}}(T))+\frac{1}{2}\Big(\nabla_2^2\mathcal L(\bar y(0),\bar y(T),\hat\ell)(X_{v,V}(T),X_{v,V}(T))
\\&&
+2\nabla_1\nabla_2\mathcal L(\bar y(0),\bar y(T);\ell)(X_{v,V}(T),V)
+\nabla_1^2\mathcal L(\bar y(0),\bar y(T);\ell)(V,V)\Big)
\\&=&\int_0^T\Big( p^\ell(t)(Y_{\sigma W}^{X_{v,V}}(t))\Big)^\prime dt+I
\\&=&\int_0^T\nabla_uH[t,\ell](\sigma(t))dt+\frac{1}{2}\int_0^T\Big(\nabla_x^2H[t,\ell](X_{v,V}(t),X_{v,V}(t))
\\&&+2\nabla_u\nabla_x H[t,\ell](X_{v,V}(t),v(t))+\nabla_u^2H[t,\ell](v(t),v(t))
\\&&-R(\tilde{p}^\ell(t),X_{v,V}(t),f[t],X_{v,V}(t))\Big)dt+I,
\end{array}$$
where 
$$\begin{array}{ll}
I=&\frac{1}{2}\nabla_1^2\mathcal L(\bar y(0),\bar y(T);\ell)(V,V)+\nabla_1\nabla_2\mathcal L(\bar y(0),\bar y(T);\ell)(X_{v,V}(T),V)
\\&+\frac{1}{2}\nabla_2^2\mathcal L(\bar y(0),\bar y(T);\ell)(X_{v,V}(T),X_{v,V}(T)),\end{array}$$
and thus (\ref{c244}) follows. 
$\Box$

%
%
%

%
%
%
%
%
%
%

\providecommand{\bysame}{\leavevmode\hbox to3em{\hrulefill}\thinspace}
\providecommand{\MR}{\relax\ifhmode\unskip\space\fi MR }
\providecommand{\MRhref}[2]{%
  \href{http://www.ams.org/mathscinet-getitem?mr=#1}{#2}
}
\providecommand{\href}[2]{#2}

\end{document}